\documentclass[final]{siamltex}

\newcommand{\uu}{{\underline u}}
\newcommand{\up}{{\underline p}}
\newcommand{\uq}{{\underline q}}
\newcommand{\uv}{{\underline v}}
\newcommand{\uw}{{\underline w}}
\newcommand{\n}{{\bf n}}
\newcommand{\PP}{{\bf P}}
\newcommand{\ddiv}{{\hbox{div}}}
\newcommand{\R}{{\mathbb R}}
\newcommand{\T}{{\mathcal T}}

\usepackage{float}
\usepackage{amsmath}
\usepackage{amssymb}
\usepackage{graphicx} 
\usepackage{psfrag} 
\usepackage{enumerate} 

\def \ds{\displaystyle}
\def \be{\begin{equation}}
\def \ee{\end{equation}}
\def \bea{\begin{eqnarray}}
\def \eea{\end{eqnarray}}
\def \bean{\begin{eqnarray*}}
\def \eean{\end{eqnarray*}}
\newtheorem{theo}{Theorem}
\newtheorem{lem}{Lemma}
\newtheorem{rem}{Remark}

\graphicspath{
{./}
{./Figures/}
{./Simulations/Example/}
{./Simulations/Test_2x1_error_alpha/}
{./Simulations/Test_mesh1234/}
{./Simulations/Test_2x2_error_alpha/}
}

\title{A new cement to glue non-conforming grids with Robin interface
conditions: 
the finite element case}

\author{Caroline Japhet\thanks{Laboratoire d'Analyse, G\'eom\'etrie et
Applications, Universit\'e Paris XIII, Avenue J-B Cl\'ement,
93430 Villetaneuse, France, E-mail: japhet@math.univ-paris13.fr} \and Yvon
Maday\thanks{Universit\'e Pierre et Marie Curie-Paris6, UMR 7598 Laboratoire
Jacques-Louis Lions, B.C. 187, Paris, F-75005 France, E-mail:
maday@ann.jussieu.fr; and Division of Applied Maths, Brown University} \and 
Fr\'ed\'eric Nataf\thanks{Universit\'e Pierre et Marie Curie-Paris6, UMR 7598 Laboratoire
Jacques-Louis Lions, B.C. 187, Paris, F-75005 France, E-mail:
nataf@ann.jussieu.fr}}

\begin{document}

\maketitle

\begin{abstract}
We design and analyze a new non-conforming domain decomposition
method based on Schwarz type approaches that allows for the use of Robin
interface conditions on non-conforming grids. The method is proven to be
well posed, and the iterative solver to converge. The error analysis is
performed in 2D  piecewise polynomials of low and high order and extended
in 3D for $P_1$ elements. Numerical results in 2D illustrate
the new method.
\end{abstract}


\section{Introduction}
Our goal in writing this paper is to propose and analyze a non-conforming
domain
decomposition generalization to P.L.~Lions initial idea \cite{Lions} in view
of an extension of the approach to optimized interface
conditions algorithms. This type of algorithm has proven indeed to be an 
efficient approach to
domain decomposition methods in the case of conforming approximations
\cite{Despres.2}, \cite{Japhet2}.
This paper presents the basic material related to so called 
zero$^{\footnotesize\mbox{th}}$
order (non optimized) method in case of finite element discretizations,
see \cite{GJMN} for a short presentation. In the companion paper
\cite{AJMN}, the case of the finite
volume discretization has been introduced and analyzed.\\
\par

We first consider the problem at the continuous level:\ \ \ Find $u$
such that
\begin{eqnarray}
\label{eq:pbgen}
{\cal L}(u)=f \hbox{ in }\Omega\\
\label{eq:pbgen2}
{\cal C}(u)=g  \hbox{ on }\partial\Omega
\end{eqnarray}
where ${\cal L}$ and ${\cal C}$ are partial differential equations.
The original Schwarz
algorithm is based on a decomposition of the domain $\Omega$ into
overlapping subdomains
and the resolution of Dirichlet boundary value problems in the 
subdomains. It has been
proposed in \cite{Lions} to use more general boundary conditions for the
problems on the subdomains in
order  to use a non-overlapping
decomposition of the domain. The convergence rate is also dramatically
increased.\\

More precisely, let $\Omega$ be a ${\cal C}^{1,1}$ (or convex polygon in 2D 
or polyhedron in 3D) domain of
$I\!\!R^d$, $d=2$ or $3$; we assume it
is decomposed into
$K$ non-overlapping
subdomains:
\begin{eqnarray}\label{de_k10}
\overline \Omega = \cup_{k=1}^{K} \overline\Omega^k.
\end{eqnarray}
We suppose that the subdomains 
$\Omega^k, \ 1 \le k \le K$ are either ${\cal C}^{1,1}$ or
convex polygons in 2D or polyhedrons in 3D. We assume also that this
decomposition is geometrically conforming in the sense
that the intersection of the closure of two different subdomains, if not
empty, is either a common vertex, a
common edge, or a common face in 3D. Let
$\n_k$ be the outward normal from
$\Omega^k$. Let $({\cal B}_{k,\ell})_{1\le k,\ell \le K, k\not= \ell}$ be
the chosen transmission conditions
on the
interface
between the subdomains (e.g. ${\cal B}_{k,\ell}={\partial\ \over
\partial \n_k}+\alpha_k$).  What we
shall call here a Schwarz type method for the problem
(\ref{eq:pbgen})-(\ref{eq:pbgen2}) is its
reformulation:
\ \ \ Find $(u_k)_{1\le k\le K}$ such that
\begin{eqnarray}
{\cal L}(u_k)=f \hbox{ in }\Omega^k \nonumber\\
{\cal C}(u_k)=g  \hbox{ on } \partial\Omega^k \cap \partial\Omega \nonumber\\
{\cal B}_{k,\ell}(u_k)={\cal B}_{k,\ell}(u_\ell)
\hbox{ on }\partial\Omega^k \cap \partial\Omega^\ell. \nonumber
\end{eqnarray}

Let us focus first on the interface conditions ${\cal B}_{k,\ell}$.
The convergence rate of associated Schwarz-type domain decomposition
methods  is very
sensitive to the choice of these transmission conditions. The use
of exact artificial (also called absorbing) boundary
conditions as interface conditions leads to an optimal number of iterations,
see \cite{Hagstrom}, \cite{Nataf.4}, \cite{GHAD2} and \cite{GHAD1}.
Indeed, for a domain decomposed into
$K$  strips, the number of iterations is $K$, see \cite{Nataf.4}. Let us
remark that this result is rather
surprising since exact absorbing conditions refer usually to truncation of
infinite domains rather than
interface conditions in domain decomposition. Nevertheless, this approach
has some drawbacks:
\begin{enumerate}
\item the explicit form of these boundary conditions is known only for constant
coefficient operators and simple geometries,
\item these boundary conditions are pseudo-differential. The cost per
iteration is
high since the corresponding discretization matrix is not sparse for the
unknowns on the boundaries of the subdomains.
\end{enumerate}
For this reason, it is usually preferred to use partial differential
approximations to the exact absorbing boundary conditions. This
approximation problem is
  classical in the field of computation on unbounded domains since the
seminal paper of
Engquist and Majda \cite{EngMajda}. The approximations correspond to ``low
frequency''
approximations of the exact absorbing boundary conditions. In domain
decomposition
methods, many authors have used them for wave propagation problems
\cite{Despres.3},
\cite{Despres.4}, \cite{Douglas}, \cite{BenDespres}, \cite{Windlund},
\cite{Keyes} and
\cite{Bourdonnaye} and in fluid dynamics
\cite{Nataf}, \cite{Quarteroni}.
Instead of using "low frequency" in space approximations to the exact
absorbing boundary
conditions, it has been proposed to design approximations which minimize
the convergence
rate of the algorithm. Such optimization of the transmission conditions
for the performance of the algorithm was done in, 
\cite{Japhet2}, \cite{JNR} for a convection-diffusion equation, where coefficients
in  second order transmission conditions where optimized. These
approximations are quite different from the "low frequency" approximations
and increase
dramatically the convergence rate of the method.

When the grids are conforming, the implementation of such interface
conditions on the
discretized problem is not too difficult. On the other hand, using
non-conforming grids is
very appealing since their use allows for parallel generation of meshes,
for local adaptive meshes and fast and
independent  solvers.
The mortar element method, first introduced in \cite{BMP}, enables the use of
non-conforming grids.
It is also well suited to the use of the so-called "Dirichlet-Neumann"
(\cite{Quarteroni}) or
"Neumann-Neumann" preconditioned conjugate gradient method applied to
the Schur
complement matrix  (\cite{Lacour}, \cite{AMW}, \cite{TosWid}). 
In the context of finite volume discretizations, it was proposed in \cite{Saas} to use a mortar type method with arbitrary interface conditions. To our knowledge, such an approach has not been extended to a finite element discretization. Moreover, the approach we present here is different and simpler. 

Our final project is to use interface conditions such as OO2
interface conditions (see \cite{Japhet2}, \cite{JNR}), this
will be
developed in a future paper.
Here we consider only interface conditions of order 0 :
${\cal B}_{k,\ell}={\partial\ \over \partial \n_k}+\alpha_k$. The approach we propose and study was introduced in \cite{GJMN} and independently implemented in \cite{Vouvakis} for the Maxwell equations but without numerical analysis.

\section{Definition of the method and of the iterative
solver.}
We consider the following problem : Find $u$ such that
\bea
\label{initial_BVP1}
(Id - \Delta)u &=& f \quad \mbox{in } \Omega \\
\label{initial_BVP2}
              u &=& 0 \quad \mbox{on } \partial{\Omega},
\eea
%
where $f$ is given in $L^2(\Omega)$. \\
The variational statement of the problem
(\ref{initial_BVP1})-(\ref{initial_BVP2})
consists in writing the problem as follows : Find $u \in H^1_0(\Omega)$
such that
\bea\label{initial_VF}
\int_{\Omega} \left(\nabla u \nabla v +uv \right) dx = \int_{\Omega} fvdx,
\quad \forall v \in H^1_0(\Omega).
\eea
%
Making use of the  domain decomposition (\ref{de_k10}), the problem
(\ref{initial_VF})
can be written as follows : Find
$u \in H^1_0(\Omega)$ such that
\bea
\sum_{k=1}^{K} \int_{\Omega^k} \left(\nabla (u_{|\Omega^k}) \
\nabla (v_{|\Omega^k}) +u_{|\Omega^k}v_{|\Omega^k} \right) dx
= \sum_{k=1}^{K} \int_{\Omega^k} f_{|\Omega^k}v_{|\Omega^k}dx,
\quad \forall v \in H^1_0(\Omega).  \nonumber
\eea
%
Let us introduce the space $H^1_*(\Omega^k)$ defined by
$$H^1_*(\Omega^k) = \{\varphi\in H^1(\Omega^k),\quad
\varphi = 0 \hbox{ over } \partial\Omega\cap\partial
\Omega^k\}.$$
It is standard to note that the space $H^1_0(\Omega)$ can then be
identified with the
subspace of the $K$-tuple $\uv=(v_1,...,v_K)$ that are continuous on the
interfaces:
\bea
V = \{\uv=(v_1,...,v_K) \in \prod_{k=1}^K H^1_*(\Omega^k), \nonumber\\ \quad
\forall k,\ell, k \ne \ell, \ 1 \le k,\ell \le K, \ v_k = v_{\ell}
\mbox{ over }
\partial{\Omega^k} \cap \partial{\Omega^{\ell}}\}. \nonumber
\eea
%
This leads to introduce also the notation of the interfaces of two adjacent
subdomains
$$\Gamma^{k,\ell} =
\partial\Omega^k\cap\partial\Omega^\ell.$$
In what follows, for the sake
of simplicity, the only fact to refer to a pair
$(k,\ell)$ preassumes that
$\Gamma^{k,\ell}$ is not empty.
The problem (\ref{initial_VF}) is then equivalent to the following one :
Find $\uu \in V$ such that
\bea
\sum_{k=1}^{K} \int_{\Omega^k} \left( \nabla u_k\nabla v_k
+u_k v_k \right) dx
= \sum_{k=1}^{K} \int_{\Omega^k} f_kv_kdx,
\quad \forall \uv \in V.  \nonumber
\eea
%
\begin{lem}\label{lemma:constraint}
For $\uv \in \prod_{k=1}^K H^1_*(\Omega^k)$, the constraint $v_k = v_{\ell}$ across the interface $\Gamma^{k,\ell}$
  is equivalent to
\bea\label{constraint}
\forall \up \equiv (p_k) \in
\prod_{k=1}^K H^{-1/2}(\partial\Omega^k) \mbox{ with }
p_k=-p_{\ell} \mbox{ over } \Gamma^{k,\ell}, \nonumber \\
\sum_{k =1}^{ K}  \
_{H^{-1/2}(\partial\Omega^k)}<p_k,v_k>_{H^{1/2}(\partial\Omega^k)}
&=& 0. 
\eea
\end{lem}
{\bf Proof} The proof is similar to the one of proposition III.1.1 in
\cite{BF91} but can't be directly derived
from this proposition. Let $ \up \equiv (p_k) \in
\prod_{k=1}^K H^{-1/2}(\partial\Omega^k)$ with
$p_k=-p_{\ell} \mbox{ over } \Gamma^{k,\ell}$,
in $(H_{00}^{1/2}(\Gamma^{k,\ell}))^{\prime}$ sense. 
Then, there exists over each $\Omega^k$ a lifting of the normal trace $p_k$ in
$H(\ddiv, \Omega^k)$. The global function $\PP$, which restriction 
to each $\Omega^k$ is defined as being equal to the lifting, belongs to
$H(\ddiv, \Omega)$ and is such that $(\PP.\n)_{|\partial \Omega^k} = p_k$.
Let now $\uv \in V$. From the previously quoted identification, we know that
there exists ${\bf v} \in H^1_0(\Omega)$ such that
${\bf v}_{|\Omega^k}=v_k$. In addition,
$$
\int_{\Omega} {\bf v} div\PP - \int_{\Omega} \PP \nabla {\bf v}  =0.
$$
On the other hand,
\bea
\int_{\Omega}  {\bf v} div\PP- \int_{\Omega} \PP \nabla {\bf v}  &=&
\sum_{k=1}^K (\int_{\Omega^k}  {\bf v} div\PP- \int_{\Omega^k} \PP \nabla {\bf v} )\nonumber\\
&=&\sum_{k=1}^K \int_{\partial \Omega^k} (\PP.\n) \uv=\sum_{k=1}^K
\int_{\partial \Omega^k} p_k v_k,\nonumber
\eea
so that (\ref{constraint}) is satisfied. \\
Reciprocally, let 
$\uv=(v_1,...,v_K) \in \prod_{k=1}^K H^1_*(\Omega^k)$ such that
(\ref{constraint}) is satisfied. Let $x \in \Gamma^{k,\ell}$, and
let $\gamma_x \subset \bar{\gamma}_x \subset \Gamma_x \subset \bar{\Gamma}_x
\subset \Gamma^{k,\ell}$ be open sets. There exists a function
$\varphi$ in ${\cal D}(\Gamma_x)$ such that $\varphi(y)=1$ for all $y$ in $\gamma_x$.
With any $q \in (H_{00}^{1/2}(\Gamma_x))^{\prime}$, let us associate
$\up \equiv (p_k)$ defined by
\bea
_{H^{-1/2}(\partial\Omega^k)}<p_k,w_k>_{H^{1/2}(\partial\Omega^k)}
&=&_{(H_{00}^{1/2}(\Gamma_x))^{\prime}}<q,\varphi w_k>_{H_{00}^{1/2}(\Gamma_x)},\forall w_k \in H^{1/2}(\partial\Omega^k), \nonumber\\
_{H^{-1/2}(\partial\Omega^{\ell})}<p_{\ell},w_{\ell}>_{H^{1/2}(\partial\Omega^{\ell})}
&=&-_{(H_{00}^{1/2}(\Gamma_x))^{\prime}}<q,\varphi w_{\ell}>_{H_{00}^{1/2}(\Gamma_x)}, \forall w_{\ell} \in H^{1/2}(\partial\Omega^{\ell}), \nonumber\\
\mbox{and } p_j&=&0, \ \forall j \neq k,\ell. \nonumber
\eea
Then, by construction, $\up \in \prod_{k=1}^K H^{-1/2}(\partial\Omega^k)$ and
$p_k=-p_{\ell} \mbox{ over } \Gamma^{k,\ell}$. Hence from
(\ref{constraint}),
$$
\sum_{k =1}^{ K}  \
_{H^{-1/2}(\partial\Omega^k)}<p_k,v_k>_{H^{1/2}(\partial\Omega^k)}
= 0. 
$$
We derive
$$
_{H^{-1/2}(\partial\Omega^k)}<p_k,v_k>_{H^{1/2}(\partial\Omega^k)}
=-_{H^{-1/2}(\partial\Omega^{\ell})}<p_{\ell},v_{\ell}>_{H^{1/2}(\partial\Omega^{\ell})},
$$
thus,
$$
_{(H_{00}^{1/2}(\Gamma_x))^{\prime}}<q,\varphi v_k>_{H_{00}^{1/2}(\Gamma_x)}
=
_{(H_{00}^{1/2}(\Gamma_x))^{\prime}}<q,\varphi v_{\ell}>_{H_{00}^{1/2}(\Gamma_x)},$$
and this is true for any $q \in (H_{00}^{1/2}(\Gamma_x))^{\prime}$,
hence
$ \varphi v_k=\varphi v_{\ell}$ over $\Gamma_x$, and thus
$$
v_k=v_{\ell} \ \mbox{over} \ \gamma_x, \quad \forall x \in \Gamma^{k,\ell}.
$$
We derive $v_k=v_{\ell}$ a.e. over $\Gamma^{k,\ell}$, which ends the proof of lemma
\ref{lemma:constraint}.\\\\
The constrained space is then defined as follows
\bea
{\cal V} = \ds\lbrace(\uv,\uq)\in \left(\prod_{k=1}^K
H^1_*(\Omega^k)\right)\times \left(\prod_{k=1}^K
H^{-1/2}(\partial\Omega^k)\right),  \nonumber\\
\quad
v_k=v_\ell\hbox{ and }q_k = - q_\ell
\hbox{ over }\Gamma^{k,\ell}\rbrace \nonumber
\eea
and problem (\ref{initial_VF}) is equivalent to the following one :
Find $(\uu,\up) \in {\cal V}$ such that\\\\
$\forall \uv \in \prod_{k=1}^KH^1_*(\Omega^k),$
\bea
\sum_{k=1}^{K} \int_{\Omega^k} \left( \nabla u_k\nabla v_k +u_kv_k \right) dx
- \sum_{k=1}^{ K}  \
_{H^{-1/2}(\partial\Omega^k)}<p_k,v_k>_{H^{1/2}(\partial\Omega^k)}
\nonumber\\
= \sum_{k=1}^{K} \int_{\Omega^k} f_kv_kdx.  \nonumber
\eea
%
Being equivalent with the original problem, where $p_k = {\partial
u\over\partial
\n_k}$ over $\partial\Omega^k$ (remind that $f$ is assumed to be in
$L^2(\Omega)$ so that
${\partial u\over\partial \n_k}$ actually belongs to
$H^{-1/2}(\partial\Omega^k)$),  this problem is
naturally well posed. This can also be directly derived from the proof
of an inf-sup condition that follows from the arguments
developed hereafter for the analysis of the iterative procedure. First, let us
describe this algorithm in the continuous case,  and then in the non
conforming discrete
case. In both  cases, we prove the convergence of the algorithm
towards the
solution of the problem.
\subsection{Continuous case}
Let us consider
the interface conditions of order 0 :
\bea
p_k + \alpha u_k = - p_{\ell} + \alpha u_{\ell} \quad \mbox{ over }
\Gamma^{k,\ell} \nonumber
\eea
%
where $\alpha$ is a given positive  real number.\\
We introduce the following algorithm : let $(u_k^n,p_k^n) \in
H^1_*(\Omega^k) \times H^{-1/2}(\partial\Omega^k)$
be an approximation of $(u,p)$ in $\Omega^k$ at step $n$.
Then, $(u_k^{n+1},p_k^{n+1})$ is the solution in
$H^1_*(\Omega^k) \times H^{-1/2}(\partial\Omega^k)$ of
\bea
\label{algo_continu}
\int_{\Omega^k} \left( \nabla u_k^{n+1}\nabla v_k
+u_k^{n+1}v_k \right) dx
- _{H^{-1/2}(\partial\Omega^k)}<p_k^{n+1},v_k>_{H^{1/2}(\partial\Omega^k)}
\nonumber\\
= \int_{\Omega^k} f_kv_kdx, \quad \forall v_k \in H^1_*(\Omega^k)  \\
\label{CI_continu}
<p_k^{n+1}+ \alpha u_k^{n+1},v_k>_{\Gamma^{k,\ell}}=
<- p_{\ell}^{n} + \alpha u_{\ell}^{n},v_k>_{\Gamma^{k,\ell}},
\quad \forall v_k \in H_{00}^{1/2}(\Gamma^{k,\ell})
\eea
%
It is obvious to remark that this series of equations results in uncoupled
problems
set on every $\Omega^k$. Recalling that $f\in L^2(\Omega)$, the strong
formulation is indeed that
\bea
-\Delta u_k^{n+1} + u_k^{n+1} = f_k\quad\hbox{ over }\Omega^k \nonumber\\
\ds{\partial u_k^{n+1}\over\partial \n_k} + \alpha u_k^{n+1} = -p_\ell^n+\alpha
u_\ell^n\quad\hbox{ over } \Gamma^{k,\ell} \nonumber\\
\label{flux_fort}
\ds p_k^{n+1} = {\partial u_k^{n+1}\over\partial \n_k} \quad\hbox{ over }
\partial\Omega^k
\eea
From this strong formulation it is
straightforward to derive by induction that if each
$p^0_k, \ k=1,...,K$, is chosen in $\prod_\ell H^{1/2}(\Gamma^{k,\ell})$,
then,
for each $k$, $1\le k\le K$, and $n\ge 0$ the solution
$u_k^{n+1}$ belongs to
$H^1(\Omega^k)$ and $p_k^{n+1}$
belongs to $\prod_\ell
H^{1/2}(\Gamma^{k,\ell})$ by standard trace results ($p_k^{n+1} = -p_\ell^n+\alpha(u^n_\ell-u_k^{n+1})$). This regularity
assumption on $p^0_k$ will be done hereafter.

  We can prove now that the algorithm (\ref{algo_continu})-(\ref{CI_continu})
converges for all $f \in L^2(\Omega)$. As the equations are linear,
we can take $f=0$. We prove the convergence in the sense that, in this
case, the
associated sequence $((u_k^n,p_k^n))_n$ satisfies
\bea
\lim_{n \longrightarrow \infty} \left( \|u_k^n\|_{H^1(\Omega^k)}
+  \|p_k^n\|_{H^{-1/2}(\partial\Omega^k)} \right)=0,
\mbox{  for } 1\le k\le K. \nonumber
\eea
%
We proceed as in (\cite{Lions},\cite{Despres.2}) by using an energy estimate that we
derive by taking
$v_k=u_k^{n+1}$ in (\ref{algo_continu}) and the use of the regularity
property that
$p_k^{n+1} \in L^2(\partial\Omega^k)$
\bea
\int_{\Omega^k} \left( |\nabla u_k^{n+1}|^2
+|u_k^{n+1}|^2 \right) dx
=  \int_{\partial\Omega^k} p_k^{n+1}u_k^{n+1} ds \nonumber
\eea
%
that can also be written
\bea
\int_{\Omega^k} \left( |\nabla u_k^{n+1}|^2
+|u_k^{n+1}|^2 \right) dx \hspace{5.6cm}  \nonumber\\
\hspace{3cm} =\sum_{\ell} \frac{1}{4\alpha} \int_{\Gamma^{k,\ell}}\left(
( p_k^{n+1}+\alpha u_k^{n+1})^2 - ( p_k^{n+1}-\alpha u_k^{n+1})^2\right)ds
  \nonumber
\eea
%
By using the interface conditions (\ref{CI_continu}) we obtain
\bea\label{estim_en}
\int_{\Omega^k} \left( |\nabla u_k^{n+1}|^2
+|u_k^{n+1}|^2 \right) dx
+\frac{1}{4\alpha}\sum_{\ell}\int_{\Gamma^{k,\ell}}
( p_k^{n+1}-\alpha u_k^{n+1})^2ds \nonumber\\
= \frac{1}{4\alpha}\sum_{\ell}\int_{\Gamma^{k,\ell}}
( - p_{\ell}^{n}+\alpha u_{\ell}^{n})^2ds
\eea
%
Let us now introduce two new quantities defined at each step $n$ :
\bea
E^n=\sum_{k=1}^K \int_{\Omega^k} \left( |\nabla u_k^{n}|^2
+|u_k^{n}|^2 \right), \nonumber
\eea
%
and
\bea
B^n = \frac{1}{4\alpha}\sum_{k=1}^K\sum_{\ell \ne k} \int_{\Gamma^{k,\ell}}
( p_k^{n}-\alpha u_k^{n})^2ds. \nonumber
\eea
%
By summing up the estimates (\ref{estim_en}) over $k=1,...,K$, we have,
\bea
E^{n+1} + B^{n+1} \le B^n, \nonumber
\eea
%
so that, by summing up these inequalities, now over $n$, we obtain :
\bea
\sum_{n=1}^{\infty} E^{n} \le B^0. \nonumber
\eea
%
We thus have $\lim_{n \longrightarrow \infty} E^n =0$.
Relation (\ref{flux_fort}) then implies :
\bea
\lim_{n \longrightarrow \infty} \|p_k^n\|_{H^{-1/2}(\partial\Omega^k)}=0,
\mbox{  for } k=1,...,K. \nonumber
\eea
which ends the proof of the convergence of the continuous algorithm.
%
\begin{theo}
Assume that $f$ is in $L^2(\Omega)$ and $(p^0_k)_{1 \le k \le K}
\in \prod_\ell H^{1/2}(\Gamma^{k,\ell})$. Then, the algorithm
(\ref{algo_continu})-(\ref{CI_continu}) converges
in the sense that
\bea
\lim_{n \longrightarrow \infty} \left( \|u_k^n - u_k\|_{H^1(\Omega^k)}
+  \|p_k^n-p_k\|_{H^{-1/2}(\partial\Omega^k)} \right)=0,
\mbox{  for } 1\le k\le K, \nonumber
\eea
%
where $u_k$ is the restriction to $\Omega^k$ of the solution $u$
to (\ref{initial_BVP1})-(\ref{initial_BVP2}), and $p_k = {\partial
u_k \over\partial \n_k}$ over $\partial\Omega^k, \ 1 \le k \le K$.
\end{theo}
\subsection{Discrete case}\label{sec.discretecase}
We  introduce now the discrete spaces. Each $\Omega^k$ is provided with
its own mesh ${\cal T}_h^k, \ 1 \le k \le K$, such that
\bea
\overline \Omega^k=\cup_{T \in {\cal T}_h^k} T. \nonumber
\eea
For $T \in {\cal T}_h^k$, let $h_T$ be the diameter of $T$
($h_T=\sup_{x,y \in T} d(x,y)$) and $h$ the discretization parameter
\bea
h=\max_{1 \le k \le K} h_k, \quad \mbox{with} \quad
h_k=\max_{T \in {\cal T}_h^k} h_T. \nonumber
\eea
At the price of (even) more techniques and care in the forthcomming analysis,
possible large variations in the norms of the solution $u_{|\Omega^k}$
can be compensated by tuning of $h_k$. 
This requires that the uniform $h$ is not used but all the analysis
is performed with $h_k$. For the sake of readability we prefer to use
$h$ instead of $h_k$.
Let $\rho_T$ be the diameter of the  circle (in 2D) or sphere (in 3D)
inscribed in $T$, then $\sigma_T=\frac{h_T}{\rho_T}$ is a measure of the
nondegeneracy of $T$. We suppose that ${\cal T}_h^k$ is uniformly regular:
there exists $\sigma$ and $\tau$ independent of $h$ such that
\bea
\forall T \in {\cal T}_h^k, \quad \sigma_T \le \sigma
\quad \mbox{and} \quad  \tau h \le h_T .\nonumber
\eea
We consider that the sets belonging to the meshes are of simplicial type
(triangles or tetrahedron), but
the analysis made hereafter can be applied as well for quadrangular or
hexahedral meshes.
Let ${\cal P}_M(T)$ denote the space of all polynomials defined over T
of total degree less than or equal to $M$. 
The finite elements are of lagrangian type, of class ${\cal C}^0$.
Then, we define over each subdomain two conforming spaces $Y_h^k$ and
$X_h^k$ by :
\bea
Y_h^k&=&\{v_{h,k} \in {\cal C}^0(\overline \Omega^k),
\ \  {v_{h,k}}_{|T} \in {\cal P}_M(T), \ \forall T \in {\cal T}_h^k \},
\nonumber\\
X_h^k&=&\{v_{h,k} \in Y_h^k, \  {v_{h,k}}_{|\partial \Omega^k \cap \partial
\Omega}=0\}.\nonumber
\eea
In what follows we assume that the mesh is designed by taking into account
the geometry of the $\Gamma^{k,\ell}$ in the sense that, the 
space of traces over each
$\Gamma^{k,\ell}$ of elements of $Y_h^k$ is a finite element space 
denoted by ${\cal Y}_h^{k,\ell}$.
Again, at the price of more notations, this assumption can be relaxed.
Let $k$ be given, the space
${\cal Y}_h^k$ is then the product space of the ${\cal Y}_h^{k,\ell}$
over each $\ell$ such that
$\Gamma^{k,\ell}\not=\emptyset$. With each such interface we associate a
subspace
$\tilde W_h^{k,\ell}$ of ${\cal Y}_h^{k,\ell}$  in the
same spirit
as in the mortar element method \cite{BMP} in 2D or \cite{BBM} and 
\cite{BraessDahmen} in 3D. 
To be more specific, let us recall the situation in 2D. If the space
$X_h^k$ consist
of continuous piecewise polynomials of degree $\le M$, then it is readily
noticed that the
restriction of $X_h^k$ to $\Gamma^{k,\ell}$ consists in finite element
functions
adapted to the (possibly curved) side $\Gamma^{k,\ell}$ of piecewise
polynomials of degree $\le
M$. This side has two end points that we denote  as
$x_0^{k,\ell}$ and $x_n^{k,\ell}$ that belong to the set of  vertices
of the
corresponding triangulation of $\Gamma^{k,\ell}$ : $x_0^{k,\ell},
x_1^{k,\ell},...,x_{n-1}^{k,\ell}, x_n^{k,\ell}$. The space $\tilde
W_h^{k,\ell}$ is
then the subspace of   those elements of ${\cal Y}_h^{k,\ell}$
that are polynomials of degree $\le M-1$
over both $[x_0^{k,\ell},
x_1^{k,\ell}]$ and $[x_{n-1}^{k,\ell}, x_n^{k,\ell}]$.
As before, the space $\tilde W_h^{k}$ is the product space of the $\tilde
W_h^{k,\ell}$ over each $\ell$ such that
$\Gamma^{k,\ell}\not=\emptyset$.
\\
The discrete constrained
space is then defined as
\begin{eqnarray}
\label{disc.const}
{\cal V}_h =  \ds\lbrace(\uu_h,\up_h)\in
\left(\prod_{k=1}^K X_h^k\right)\times
\left(\prod_{k=1}^K \tilde W_h^{k}\right), \nonumber\\
\ \qquad\int_{\Gamma^{k,\ell}}((p_{h,k}+\alpha u_{h,k})-(-p_{h,\ell}+\alpha
u_{h,\ell})
)\psi_{h,k,\ell}
= 0,\
\forall \psi_{h,k,\ell} \in \tilde W_h^{k,\ell}
\rbrace,
\end{eqnarray}
and the discrete problem is the following one :
Find $(\uu_h,\up_h) \in {\cal V}_h$ such that\\\\
$\forall \uv_h=(v_{h,1},...v_{h,K}) \in \prod_{k=1}^K X_h^k,$
\bea\label{pbdiscret}
\sum_{k=1}^{K} \int_{\Omega^k} \left( \nabla u_{h,k}\nabla v_{h,k} +u_{h,k}
v_{h,k} \right) dx
- \sum_{k=1}^{K} \int_{\partial\Omega^k} p_{h,k} v_{h,k} ds
= \sum_{k=1}^{K} \int_{\Omega^k} f_k v_{h,k} dx.\hspace{10mm}
\eea
%
%
We introduce the discrete algorithm : let $(u_{h,k}^n,p_{h,k}^n) \in
X_h^k \times \tilde W_h^{k}$
be a discrete approximation of $(u,p)$ in $\Omega^k$ at step $n$.
Then, $(u_{h,k}^{n+1},p_{h,k}^{n+1})$ is the solution in $X_h^k
\times\tilde W_h^k$ of
\bea
\label{algo_discret}
\int_{\Omega^k} \left( \nabla u_{h,k}^{n+1}\nabla v_{h,k}
+u_{h,k}^{n+1}v_{h,k} \right) dx -
  \int_{\partial\Omega^k}p_{h,k}^{n+1} v_{h,k} ds
= \int_{\Omega^k} f_kv_{h,k}dx  ,\ \forall v_{h,k}\in X_h^k \hspace{10mm}\\
\label{CI_discret}
\int_{\Gamma^{k,\ell}} (p_{h,k}^{n+1}+ \alpha u_{h,k}^{n+1})\psi_{h,k,\ell} =
\int_{\Gamma^{k,\ell}} ( -p_{h,\ell}^{n} + \alpha
u_{h,\ell}^{n}) \psi_{h,k,\ell} ,
\quad \forall \psi_{h,k,\ell} \in  \tilde W_h^{k,\ell}.\hspace{10mm}
\eea
In order to analyze the convergence of this iterative scheme, we have to
precise the norms that can be used on the Lagrange multipliers $\up_h$.
For any $\up \in \prod_{k=1}^{K}L^2(\partial \Omega^k)$, in addition to
the natural norm, we can define two better suited norms as follows
\bea
\|\up\|_{- {1 \over 2},*} = \left(\sum_{k=1}^K \sum_{\scriptstyle \ell=1
\atop{\atop \scriptstyle \ell \ne k}}^K
\|p_k\|_{H^{-{1\over 2}}_*(\Gamma^{k,\ell})}^2 \right)^{1 \over 2},
\nonumber
\eea
where $\|.\|_{H^{-{1\over 2}}_*(\Gamma^{k,\ell})}$ stands for the dual norm of
${H^{{1\over 2}}_{00}(\Gamma^{k,\ell})}$ and
\bea
\|\up\|_{-{1\over 2}} = \left(\sum_{k=1}^K
\|p_k\|_{H^{-{1\over 2}}(\partial \Omega^k)}^2 \right)^{1 \over 2}.
\nonumber
\eea

We also need a stability result for the Lagrange multipliers, and refer to
\cite{BB} in 2D and to the annex in 3D, in which it is
proven that,

\begin{lem}\label{lem.faker}
There exists a constant
$c_*$ such that, for
any $p_{h,k,\ell}$ in $\tilde W_h^{k,\ell}$, there exists an element
$w^{h,k,\ell}$ in
$X_h^k$ that
vanishes  over $\partial\Omega^k\setminus\Gamma^{k,\ell}$ and satisfies
\bea
\label{stab1}
\int_{\Gamma^{k,\ell}} p_{h,k,\ell} w^{h,k,\ell} \ge
\|p_{h,k,\ell}\|^2_{H^{-{1\over
2}}_*(\Gamma^{k,\ell})}
\eea
with a bounded norm
\bea
\label{stab2}
\|w^{h, k,\ell}\|_{H^1(\Omega^k)} \le c_* \|p_{h,k,\ell}\|_{H^{-{1\over
2}}_*(\Gamma^{k,\ell})}.
\eea
\end{lem}

Let $\pi_{k,\ell}$ denote the orthogonal projection operator from
$L^2(\Gamma^{k,\ell})$
onto $\tilde W_h^{k,\ell}$. Then, for $v \in L^2(\Gamma^{k,\ell})$,
$\pi_{k,\ell}(v)$ is the unique element of $\tilde W_h^{k,\ell}$
such that
\begin{eqnarray}\label{eq:defpi}
\int_{\Gamma^{k,\ell}} (\pi_{k,\ell}(v)-v)\psi=0, \quad \forall \psi \in
\tilde W_h^{k,\ell}.
\end{eqnarray}

We are now in a position to prove the convergence of the iterative scheme

\begin{theo}\label{theo2}
Let us assume that $\alpha h \le c$, for some constant $c$ small enough.
Then, the discrete problem (\ref{pbdiscret}) has a unique solution
$(\uu_h,\up_h) \in {\cal V}_h$.
The algorithm (\ref{algo_discret})-(\ref{CI_discret}) is well posed and
converges
in the sense that
\bea
\lim_{n \longrightarrow \infty} \left( \|u_{h,k}^n - u_{h,k}\|_{H^1(\Omega^k)}
+ \sum_{\ell\neq k} \|p_{h,k,\ell}^n-p_{h,k,\ell}\|_{H^{-{1\over
2}}_*(\Gamma^{k,\ell})} \right)=0,
\mbox{  for } 1\le k\le K. \nonumber
\eea
%
\end{theo}
%
{\bf Proof}. For the sake of convenience, we drop out the index $h$ in what
follows.
We first assume that problems (\ref{pbdiscret}) and
(\ref{algo_discret})-(\ref{CI_discret}) are well posed and
proceed as in the continuous case and assume that $f=0$.
From (\ref{eq:defpi}) we have

$$\forall v_k \in L^2(\Gamma^{k,\ell}), \quad
\int_{\Gamma^{k,\ell}} p_k^{n+1} v_k = \int_{\Gamma^{k,\ell}} p_k^{n+1}
\pi_{k,\ell}(v_k),$$
and (\ref{CI_discret}) also reads
\bea
p_k^{n+1}+\alpha \pi_{k,\ell} (u_k^{n+1})=
\pi_{k,\ell} (-p_{\ell}^n+\alpha u_{\ell}^n) \quad \mbox{over }
\Gamma^{k,\ell}. \nonumber
\eea
%
By taking $v_k=u_k^{n+1}$ in (\ref{algo_discret}), we thus have
\bea
\int_{\Omega^k} \left( |\nabla u_k^{n+1}|^2
+|u_k^{n+1}|^2 \right) dx \hspace{7cm}\nonumber\\
=  \sum_{\ell} \frac{1}{4\alpha} \int_{\Gamma^{k,\ell}}\left(
( p_k^{n+1}+\alpha \pi_{k,\ell} (u_k^{n+1}))^2 - ( p_k^{n+1}-\alpha
\pi_{k,\ell} (u_k^{n+1}))^2\right)ds.
  \nonumber
\eea
%
Then, by using the interface conditions (\ref{CI_discret}) we obtain
\bea
\int_{\Omega^k} \left( |\nabla u_k^{n+1}|^2
+|u_k^{n+1}|^2 \right) dx
+\frac{1}{4\alpha}\sum_{\ell}\int_{\Gamma^{k,\ell}}
( p_k^{n+1}-\alpha \pi_{k,\ell}(u_k^{n+1}))^2ds\nonumber\\
= \frac{1}{4\alpha}\sum_{\ell}\int_{\Gamma^{k,\ell}}
  (\pi_{k,\ell}(p_{\ell}^{n}-\alpha u_{\ell}^{n}))^2 ds. \nonumber
\eea
It is straightforward to note that
\bea
\int_{\Gamma^{k,\ell}}
(\pi_{k,\ell}(p_{\ell}^{n}-\alpha u_{\ell}^{n}))^2ds
\le \int_{\Gamma^{k,\ell}}
( p_{\ell}^{n}-\alpha  u_{\ell}^{n})^2 ds\nonumber\\
=  \int_{\Gamma^{k,\ell}}(p_{\ell}^{n} -\alpha \pi_{\ell,k}(u_{\ell}^{n}) +
\alpha \pi_{\ell,k}(u_{\ell}^{n}) -
\alpha  u_{\ell}^{n})^2ds\nonumber\\
= \int_{\Gamma^{k,\ell}}(p_{\ell}^{n} -\alpha \pi_{\ell,k}(u_{\ell}^{n}))^2 +
\alpha ^2(\pi_{\ell,k}(u_{\ell}^{n}) - u_{\ell}^{n})^2ds \nonumber
\eea
since $(Id-\pi_{\ell,k})(u_{\ell}^{n})$ is orthogonal to any element in
$\tilde W_h^{\ell,k}$. We then recall that (see \cite{BMP} in 2D and
\cite{BBM} or \cite {BraessDahmen} equation (5.1) in 3D)
\bea
\label{eq:propr-pilk}
\int_{\Gamma^{k,\ell}}(\pi_{\ell,k}(u_{\ell}^{n}) - u_{\ell}^{n})^2ds
\le c h \|u_\ell^n\|_{H^{1/2}(\Gamma^{k,\ell})}^2\\
\le c h \| u_\ell^n\|_{H^1(\Omega^{\ell})}^2.\nonumber
\eea
With similar notations as those introduced in the
continuous case, we deduce
\bea
E^{n+1} + B^{n+1} \le c \alpha h E^n + B^n \nonumber
\eea
%
and we conclude as in the continuous case: if $c \alpha h < 1$ then
$\lim_{n\rightarrow\infty}E^n = 0$. The
convergence of $u_k^n$ towards 0 in the $H^1$ norm follows. The convergence
of $p_k^n$ in the
$H^{-{1\over
2}}_*(\Gamma^{k,\ell})$ norm is then derived from (\ref{stab1}) and
(\ref{algo_discret}).
Note that by having $f=0$ and $(u^n,p^n)=0$ prove that
$(u^{n+1},p^{n+1})=0$ from which we derive that the square problem
(\ref{algo_discret})-(\ref{CI_discret}) is uniquely solvable hence well posed.
Similarly, having $f=0$ and getting rid of the superscripts $n$ and $n+1$
in the previous proof gives (with obvious notations) :
\bea
E + B \le c \alpha h E + B.\nonumber
\eea
%
The well posedness of (\ref{pbdiscret}) then results with similar arguments.
\\\\
We shall address, in what follows this question through a  more direct
proof. This will allow in turn, to provide some analysis of the approximation
properties of this
scheme.

\section{Numerical Analysis.}

\subsection{Well posedness.}

The first step in this error analysis is to prove the stability of the
discrete problem
and thus its well posedness. Let us introduce over $(\prod_{k=1}^K
H^1_*(\Omega^k)\times \prod_{k=1}^K  L^2(\partial\Omega^k))\times
\prod_{k=1}^K  H^1_*(\Omega^k)$ the bilinear form
\bea
\label{fbs_discret}
\tilde a((\uu,\up), \uv)) = \sum_{k=1}^K
\int_{\Omega^k} \left( \nabla u_k\nabla v_k
+u_k v_k \right) dx -\sum_{k=1}^K
  \int_{\partial\Omega^k}p_k v_k ds.
\eea
The space $\prod_{k=1}^K H^1_*(\Omega^k)$ is endowed with the norm
\bea
\|\uv\|_* = \left(\sum_{k=1}^K \|v_k\|_{H^1(\Omega^k)}^2 \right)^{1 \over
2}.
\nonumber
\eea
\begin{lem}\label{lem.infsup}

There exists a constant $\beta>0$ independent of $h$ such that
\bea
\label{inf-sup_discret}
\forall (\uu_h,\up_h) \in {\cal V}_h ,\ \exists \uv_h\in\prod_{k=1}^K
X_h^k, \hspace{3cm}\nonumber\\
\tilde a((\uu_h,\up_h), \uv_h)) \ge
\beta (\|\uu_h\|_*+ \|\up_h\|_{-{1\over 2},*}) \|\uv_h\|_*.
\eea
%
Moreover, we have the continuity argument : there exists a constant $c>0$
such that
\bea
\forall (\uu_h,\up_h) \in {\cal V}_h ,\ \forall \uv_h\in\prod_{k=1}^K
X_h^k, \hspace{3cm}\nonumber\\
\label{ineq:continuity}
\tilde a((\uu_h,\up_h), \uv_h)) \le c (\| \uu_h \|_* +
\|\up_h\|_{-{1 \over 2}}) (\| \uv_h \|_*).
\eea
\end{lem}
%
{\bf Proof of lemma \ref{lem.infsup}:}
In (16) and (17), we have introduced local $H^1_0(\Gamma_{k,\ell})$
functions that can be put together in order to provide  an
element $\uw_h$ of $\prod_{k=1}^K X_h^k$ that satisfies
\bea
\sum_{k=1}^K \int_{\partial\Omega^k}p_k w_k ds \ge 
\|\up_h\|_{-{1 \over 2},*}^2. \nonumber
\eea
Let us now choose a real number $\gamma$,
$0<\gamma<{2\over c_*^2}$ (where $c_*$ is introduced in (\ref{stab2})) and
choose $\uv_h = \uu_h - \gamma\uw_h$ in
(\ref{fbs_discret})  that yields
\bea
\label{1.1}
\tilde a((\uu_h,\up_h), \uv_h)) = \sum_{k=1}^K
\int_{\Omega^k} \left( \nabla u_k \nabla (u_k -\gamma w_k)
+u_k (u_k -\gamma w_k) \right) dx \nonumber\\
  - \sum_{k=1}^K\int_{\partial\Omega^k}p_k(u_k -\gamma w_k) ds
\eea
As already noticed in section \ref{sec.discretecase}, we can write
\bea
\int_{\Gamma^{k,\ell}}p_k u_k  ds &=& {1\over 4\alpha}
\int_{\Gamma^{k,\ell}}((p_k + \alpha
\pi_{k,\ell}(u_k))^2   - (p_k - \alpha \pi_{k,\ell}(u_k))^2 ) ds \nonumber\\
&=& {1\over 4\alpha} \int_{\Gamma^{k,\ell}}((\pi_{k,\ell}( - p_\ell + \alpha
  u_\ell))^2   - (p_k - \alpha \pi_{k,\ell}(u_k))^2 ) ds \nonumber\\
&\le & {1\over 4\alpha} \int_{\Gamma^{k,\ell}}((p_\ell - \alpha
u_\ell)^2   - (p_k - \alpha \pi_{k,\ell}(u_k))^2 ) ds \nonumber\\
&\le & {1\over 4\alpha} \int_{\Gamma^{k,\ell}}((p_\ell - \alpha
\pi_{\ell,k}(u_\ell))^2   - (p_k - \alpha \pi_{k,\ell}(u_k))^2)ds
\nonumber\\
&&+ \ {1\over 4\alpha} \int_{\Gamma^{k,\ell}}
\alpha^2 (\pi_{\ell,k}(u_\ell)-u_\ell)^2 ds \quad \quad \nonumber
\eea
so that
\bea
\sum_{k=1}^K
  \int_{\partial\Omega^k}p_k u_k ds \le {\alpha \over 4} \sum_{k=1}^K
\int_{\partial\Omega^k} (u_k - \pi_{k,\ell}(u_k))^2 ds \le c \alpha h
\| \uu_h \|^2_* . \nonumber
\eea
Going back to (\ref{1.1}) yields
\bea
\tilde a((\uu_h,\up_h), \uv_h) &\ge& (1-c \alpha h) \| \uu_h \|^2_*  - \gamma
\| \uu_h
\|_*
\| \uw_h \|_* + \gamma \|\up_h\|_{-{1 \over 2},*}^2 \nonumber\\
&\ge& ({1\over 2}-c\alpha h) \| \uu_h \|^2_* + \gamma \|\up_h\|_{-{1 \over
2},*}^2 - {\gamma^2\over 2}
\| \uw_h \|_*^2
\nonumber\\
&\ge& ({1\over 2}-c\alpha h) \| \uu_h \|^2_* + (\gamma -
{\gamma^2c^2_*\over 2})
\|\up_h\|_{-{1 \over 2},*}^2 .\nonumber
\eea
Due to the choice of $\gamma$, we know that, for $\alpha h$ small enough,
(\ref{inf-sup_discret}) holds.
The continuity (\ref{ineq:continuity}) follows from standard arguments
(note that the norm on the right hand side of (\ref{ineq:continuity}) is
not the
$\|.\|_{-{1 \over 2},*}$--norm), which ends the proof of lemma
\ref{lem.infsup}.
\\\\
From this lemma, we have the following result :
\begin{theo}
Under the hypothesis of theorem \ref{theo2},
the discrete problem (\ref{pbdiscret}) has a unique solution
$(\uu_h,\up_h) \in {\cal V}_h$, and there exists a constant $c>0$
such that
\bea
\| \uu_h\|_* +  \|\up_h\|_{-{1 \over 2},*}  \le c \| f \|_{L^2(\Omega)}.
\nonumber
\eea
\end{theo}
%
\noindent\\
From lemma \ref{lem.infsup}, we are also in position to state that
for any $(\tilde\uu_h,\tilde \up_h)\in{\cal V}_h$,
\bea\label{estimuuh}
\| \uu - \uu_h\|_* +  \|\up - \up_h\|_{-{1 \over 2},*}  \le
c (\| \uu - \tilde\uu_h\|_* +  \|\up - \tilde\up_h\|_{-{1 \over 2}})
\eea
and we are naturally led to the analysis of the best fit of $(\uu,\up)$
by elements
in ${\cal V}_h$.

\subsection{Analysis of the best fit in 2D}\label{sec.bestfit2D}
In this part we analyze the best approximation of $(\uu,\up)$ by
elements in ${\cal V}_h$. As the proof is very technical for the analysis
of the best fit, we restrict ourselves in this section
to the complete analysis of the 2D situation. The extension to 3D first
order approximation is postponed to the next subsection.

The first step in the analysis is to prove the following lemma

\begin{lem}\label{lem_1}
Assume the degree of the finite element approximation $M\le 13$, there
exists two constants $c_1>0$ and $c_2>0$ independent of $h$
such that for all
$\eta_{\ell,k}$ in
${\cal Y}_h^{\ell,k}\cap H_0^1(\Gamma^{k,\ell})$, there exists an element
$\psi_{\ell,k}$ in
$\tilde W_h^{\ell,k}$,   such that
\bea\label{injectif}
\int_{\Gamma^{k,\ell}}(\eta_{\ell,k}+\pi_{k,\ell}(\eta_{\ell,k}))\psi_{\ell,k}
\ge c_1\|\eta_{\ell,k} \|_{L^2(\Gamma^{k,\ell})}^2,
\eea
%
\bea\label{stable}
\| \psi_{\ell,k} \|_{L^2(\Gamma^{k,\ell})} \le
c_2 \| \eta_{\ell,k} \|_{L^2(\Gamma^{k,\ell})}.
\eea
\end{lem}
%
Then, we can prove the following interpolation estimates :
\begin{theo}
\label{best-fit}
For any $u\in H^2(\Omega)\cap H^1_0(\Omega)$, such that $u_k=u_{|\Omega^k}\in
H^{2+m}(\Omega^k)$, $1\le k\le K$ with
$M-1\ge m \ge 0$,  $\uu=(u_k)_{1\le k\le K}$
and let
$p_{k,\ell}=\frac{\partial u}{\partial {\bf n}_k}$
over each $\Gamma^{k,\ell}$.
Then there exists an element $\tilde{\uu}_h$ in
$\prod_{k=1}^K X_h^k$
and $\tilde{\up}_h=(\tilde{p}_{k \ell h}), \  \tilde{p}_{k \ell h}
\in \tilde W_h^{k,\ell}$
such that $(\tilde{\uu}_h,\tilde{\up}_h)$ satisfy the coupling condition
(\ref{disc.const}), and
\bea
\| \tilde{\uu}_h -\uu\|_*
&\le& c h^{1+m} \sum_{k=1}^K \| u_k \|_{H^{2+m}(\Omega^k)}
  +{c h^m \over \alpha} \sum_{k < \ell} \| p_{k,\ell} \|_{H^{{1 \over
2}+m}(\Gamma^{k,\ell})}
\nonumber\\ \nonumber\\
\| \tilde{p}_{k \ell h} - p_{k,\ell} \|_{H^{-{1 \over 2}}(\Gamma^{k,\ell})}
&\le& c\alpha h^{2+m} (\|u_k\|_{H^{2+m}(\Omega^k)}
+\|u_{\ell}\|_{H^{2+m}(\Omega^{\ell})}) \nonumber\\
&& + \ c h^{1+m} \| p_{k,\ell} \|_{H^{{1 \over 2}+m}(\Gamma^{k,\ell})}
\nonumber
\eea
where $c$ is a constant independent of $h$ and $\alpha$.
\end{theo}
%
If we assume more regularity on the normal derivatives on the interfaces,
we have
\begin{theo}
\label{best-fit.2}
Let $u\in H^2(\Omega)\cap H^1_0(\Omega)$, such that $u_k=u_{|\Omega^k}\in
H^{2+m}(\Omega^k)$, $1\le k\le K$ with
$M-1\ge m \ge 0$,  $\uu=(u_k)_{1\le k\le K}$
and
$p_{k,\ell}=\frac{\partial u}{\partial {\bf n}_k}$
is in $H^{{3 \over 2}+m}(\Gamma_{k,\ell})$.
Then there exists $\tilde{\uu}_h$ in
$\prod_{k=1}^K X_h^k$
and $\tilde{\up}_h=(\tilde{p}_{k \ell h}), \  \tilde{p}_{k \ell h} \in
\tilde W_h^{k,\ell}$
such that $(\tilde{\uu}_h,\tilde{\up}_h)$ satisfy
(\ref{disc.const}), and
\bea
\| \tilde{\uu}_h -\uu\|_*
&\le& c h^{1+m} \sum_{k=1}^K \| u_k \|_{H^{2+m}(\Omega^k)} \nonumber \\
& & +{c h^{m+1} \over \alpha}(\log h)^{\beta(m)} \sum_{k < \ell} \| p_{k,\ell}
\|_{H^{{3
\over 2}+m}(\Gamma^{k,\ell})}
\nonumber\\ \nonumber\\
\| \tilde{p}_{k \ell h} - p_{k,\ell} \|_{H^{-{1 \over 2}}(\Gamma^{k,\ell})}
&\le& c\alpha h^{2+m} (\|u_k\|_{H^{2+m}(\Omega^k)}
+\|u_{\ell}\|_{H^{2+m}(\Omega^{\ell})}) \nonumber\\
&& + \ c h^{1+m} (\log h)^{\beta(m)}\| p_{k,\ell} \|_{H^{{1 \over
2}+m}(\Gamma^{k,\ell})}\nonumber
\eea
where $c$ is a constant independent of $h$ and $\alpha$, and $\beta(m)=0$ if
$m\le M-2$ and $\beta(m)=1$ if $m=M-1$.
\end{theo}
%

\noindent
\\
{\bf Proof of lemma \ref{lem_1}}:
We first give the {\bf proof for $P_1$ finite
element}.
The general proof is quite technical and is based on Lemma~\ref{lem:etapsi_base},
given in Appendix~\ref{sec:Ext2D}.
  Remind that we have denoted as
$x_0^{\ell,k},
x_1^{\ell,k},...,x_{n-1}^{\ell,k}, x_n^{\ell,k}$ the vertices of the
triangulation of $\Gamma^{\ell,k}$ that belong to $\Gamma^{\ell,k}$.
To any $\eta_{\ell,k}$ in
${\cal Y}_h^{\ell,k}\cap H_0^1(\Gamma^{k,\ell})$ we then associate the
element $\psi_{\ell,k}$ in
$\tilde W_h^{\ell,k}$ as follows
\bea\label{decomp.psi}
\psi_{\ell,k}=\left\{
\begin{array}{l}
{\eta_{\ell,k} (x_1^{\ell,k}-x_0^{\ell,k})\over(x-x_0^{\ell,k})}
\hbox{ over }
]x_0^{\ell,k}, x_{ 1}^{\ell,k}[\\
\eta_{\ell,k}\hbox{ over }
]x_1^{\ell,k}, x_{n-1}^{\ell,k}[\\
{\eta_{\ell,k}(x_n^{\ell,k}-x_{n-1}^{\ell,k})\over(x_n^{\ell,k}-x)}
\hbox{ over }
]x_{n-1}^{\ell,k}, x_{n}^{\ell,k}[\\
\end{array}\right.\nonumber
\eea
recalling that all norms are equivalent over the space of polynomials of
degree $1$ we easily deduce that
there exists a constant $c$ such that
\bea
\label{estim_psieta}
\| \psi_{\ell,k} \|_{L^2(\Gamma^{k,\ell})} \le
c \| \eta_{\ell,k} \|_{L^2(\Gamma^{k,\ell})}.\nonumber
\eea
Moreover, it is straightforward to derive
\bea
\int_{\Gamma^{k,\ell}}(\eta_{\ell,k}+\pi_{k,\ell}(\eta_{\ell,k}))\psi_{\ell,k}
&=&\int_{\Gamma^{k,\ell}} \eta_{\ell,k}\psi_{\ell,k} + \int_{\Gamma^{k,\ell}}
(\pi_{k,\ell}(\eta_{\ell,k}))^2 \nonumber\\
&&+ \int_{\Gamma^{k,\ell}}
\pi_{k,\ell}(\eta_{\ell,k})(\psi_{\ell,k} - \eta_{\ell,k}). \nonumber
\eea

Then, by using the relation
\bea
\pi_{k,\ell}(\eta_{\ell,k})(\psi_{\ell,k} - \eta_{\ell,k})
\ge -{1 \over 2} (\pi_{k,\ell}(\eta_{\ell,k}))^2-{1 \over
2}(\psi_{\ell,k} - \eta_{\ell,k})^2 \nonumber
\eea
leads to
\bea\label{int.eta1}
\int_{\Gamma^{k,\ell}}(\eta_{\ell,k}+\pi_{k,\ell}(\eta_{\ell,k}))\psi_{\ell,k}
\ge \int_{\Gamma^{k,\ell}} \eta_{\ell,k}\psi_{\ell,k}  + {1 \over 2}
\int_{\Gamma^{k,\ell}}
(\pi_{k,\ell}(\eta_{\ell,k}))^2 -{1 \over 2} \int_{\Gamma^{k,\ell}}
(\psi_{\ell,k} - \eta_{\ell,k})^2. \nonumber
\eea
We realize now that, over the first interval,
$$\int_{]x_0^{\ell,k}, x_{ 1}^{\ell,k}[} (\eta_{\ell,k}\psi_{\ell,k} -{1
\over 2}
(\psi_{\ell,k} - \eta_{\ell,k})^2 )= \int_{]x_0^{\ell,k}, x_{ 1}^{\ell,k}[}
({(x-x_0^{\ell,k})\over(x_1^{\ell,k}-x_0^{\ell,k})}-{1\over 2}
{(x-x_1^{\ell,k})^2\over(x_1^{\ell,k}-x_0^{\ell,k})^2})\psi_{\ell,k}^2 $$
noticing that $$\int_{]x_0^{\ell,k}, x_{ 1}^{\ell,k}[}({(x-x_0^{\ell,k})\over(x_1^{\ell,k}-x_0^{\ell,k})}-{1\over 2}
{(x-x_1^{\ell,k})^2\over(x_1^{\ell,k}-x_0^{\ell,k})^2})= \int_{]x_0^{\ell,k}, x_{ 1}^{\ell,k}[}{(x-x_0^{\ell,k})^2\over(x_1^{\ell,k}-x_0^{\ell,k})^2},$$
by recalling that $\psi_{\ell,k}$ is constant on
$]x_0^{\ell,k}, x_{ 1}^{\ell,k}[$ and $]x_{n-1}^{\ell,k}, x_{n}^{\ell,k}[$,
we derive that
\bea
\int_{\Gamma^{k,\ell}}(\eta_{\ell,k}+\pi_{k,\ell}(\eta_{\ell,k}))\psi_{\ell,k}
\ge \int_{\Gamma^{k,\ell}} \eta_{\ell,k}^2
\nonumber
\eea
which ends the proof of lemma \ref{lem_1}.
\\\\
{\bf PROOF of Lemma \ref{lem_1} in the general case}
Using the definition of $\pi_{k,\ell}$, (\ref{eq:defpi}), it is straightforward to derive
\bea
\int_{\Gamma^{k,\ell}}(\eta_{\ell,k}+\pi_{k,\ell}(\eta_{\ell,k}))\psi_{\ell,k}
&=&\int_{\Gamma^{k,\ell}} \eta_{\ell,k}\psi_{\ell,k} + \int_{\Gamma^{k,\ell}}
(\pi_{k,\ell}(\eta_{\ell,k}))^2 \nonumber\\
&&+ \int_{\Gamma^{k,\ell}}
\pi_{k,\ell}(\eta_{\ell,k})(\psi_{\ell,k} - \eta_{\ell,k}). \nonumber
\eea

Then, using the relation
\bea
\pi_{k,\ell}(\eta_{\ell,k})(\psi_{\ell,k} - \eta_{\ell,k})
\ge - (\pi_{k,\ell}(\eta_{\ell,k}))^2-{1 \over
4}(\psi_{\ell,k} - \eta_{\ell,k})^2 \nonumber
\eea
leads to
\bea
\int_{\Gamma^{k,\ell}}(\eta_{\ell,k}+\pi_{k,\ell}(\eta_{\ell,k}))\psi_{\ell,k}
\ge \int_{\Gamma^{k,\ell}} \eta_{\ell,k}\psi_{\ell,k}  -{1 \over 4}
\int_{\Gamma^{k,\ell}} (\psi_{\ell,k} - \eta_{\ell,k})^2. \nonumber
\eea
  Remind that we have denoted as
$x_0^{\ell,k},
x_1^{\ell,k},...,x_{n-1}^{\ell,k}, x_n^{\ell,k}$ the vertices of the
triangulation of $\Gamma^{\ell,k}$ that belong to $\Gamma^{\ell,k}$.
By Lemma~\ref{lem:etapsi_base} of appendix B, and an easy 
scaling argument,
there exists $c, C >0$, $\psi_1$ from $[x_0^{\ell,k},x_1^{\ell,k}]$ 
into $\R$ and
$\psi_n$ from
$[x_{n-1}^{\ell,k},x_n^{\ell,k}]$ into $\R$ such that
\[
\| \psi_1\|_{L^2(x_0^{\ell,k},x_1^{\ell,k})} + \|
\psi_n\|_{L^2(x_{n-1}^{\ell,k},x_n^{\ell,k})} \le C^2  (\|
\eta\|_{L^2(x_0^{\ell,k},x_1^{\ell,k})} +
\|
\eta\|_{L^2(x_{n-1}^{\ell,k},x_n^{\ell,k})}),
\]
$\psi_1(x_1^{\ell,k})=\eta(x_1^{\ell,k})$,
$\psi_N(x_{n-1}^{\ell,k})=\eta(x_{n-1}^{\ell,k})$ and
\[
\int_{x_0^{\ell,k}}^{x_1^{\ell,k}} (\eta \psi_1  -{1 \over 4}
  (\psi_1 - \eta)^2)
+ \int_{x_{n-1}^{\ell,k}}^{x_n^{\ell,k}} (\eta \psi_n  -{1 \over 4}
  (\psi_n - \eta)^2)
\ge c (\int_{x_0^{\ell,k}}^{x_1^{\ell,k}} \eta^2+
\int_{x_{n-1}^{\ell,k}}^{x_n^{\ell,k}}
\eta^2 ).
\]
Taking
  $\psi_{\ell,k}$ in
$\tilde W_h^{\ell,k}$ as follows
\bea
\psi_{\ell,k}=\left\{
\begin{array}{l}
  \psi_1
\hbox{ over }
]x_0^{\ell,k}, x_{ 1}^{\ell,k}[\\
\eta_{\ell,k}\hbox{ over }
]x_1^{\ell,k}, x_{n-1}^{\ell,k}[\\
\psi_n
\hbox{ over }
]x_{n-1}^{\ell,k}, x_{n}^{\ell,k}[\\
\end{array}\right.\nonumber
\eea
proves Lemma~\ref{lem_1} with $c_1=\min(1,c)$ and $c_2=\max(1,C)$.
\\\\
{\bf Proof of theorem \ref{best-fit}:} In order to prove this theorem, let
us build an element that will belong to the
discrete space and will be as close as the expected error to the solution.
Let $u_{kh}^1$ be the
unique element of $X_h^k$ defined as follows :
\begin{itemize}
\item
$(u_{kh}^1)_{|\partial \Omega^k}$ is the best fit of $u_k$ over $\partial \Omega^k$ in 
${\cal Y}_h^{k,\ell}$,
\item
$u_{kh}^1$ at the inner nodes of the triangulation (in $\Omega^k$) coincide
with the interpolate of $u_k$.
\end{itemize}
Then, it satisfies
\bea\label{bestfit_u.2}
\| u_{kh}^1 - u_k\|_{L^2(\partial \Omega^k)}
\le c h^{{3 \over 2}+m} \|u_k\|_{H^{2+m}(\Omega^k)},
\eea
%
from which we deduce that
\bea
\label{bestfit_u}
\| u_{kh}^1 - u_k\|_{L^2(\Omega^k)} + h\| u_{kh}^1 - u_k\|_{H^1(\Omega^k)}
\le c h^{2+m} \|u_k\|_{H^{2+m}(\Omega^k)},
\eea
%
and, from Aubin-Nitsche estimate
\bea
\label{bestfit_u.3}
\| u_{kh}^1 - u_k\|_{H^{-{1 \over 2}}(\Gamma^{k,\ell})}
\le c h^{2+m} \|u_k\|_{H^{2+m}(\Omega^k)}.
\eea
%
%
We define then separately the best fit $p_{k \ell h}^1$ of
$p_{k,\ell}=\frac{\partial u}{\partial {\bf n}_k}$
over each
$\Gamma^{k,\ell}$ in $\tilde W_h^{k,\ell}$. These elements satisfy 
for $0\le m\le M-1$
the error estimate
\bea
\label{bestfit_p}
\| p_{k \ell h}^1 - p_{k,\ell} \|_{L^2(\Gamma^{k,\ell})}
&\le& c h^{{1 \over 2}+m} \| p_{k,\ell} \|_{H^{{1 \over 2}+m}(\Gamma^{k,\ell})}
\\
\label{bestfit_p.2}
\| p_{k \ell h}^1 - p_{k,\ell} \|_{H^{-{1 \over 2}}(\Gamma^{k,\ell})}
&\le& c h^{1+m} \| p_{k,\ell} \|_{H^{{1 \over 2}+m}(\Gamma^{k,\ell})}.
\eea
But there is very few chance that $(\uu_h^1,\up_h^1)$ satisfy the coupling
condition (\ref{disc.const}). This element of
$\left(\prod_{k=1}^K X_h^k\right)\times
\left(\prod_{k=1}^K \tilde W_h^{k}\right)$
misses (\ref{disc.const}) of
elements $\epsilon_{k,\ell}$ and $\eta_{\ell,k}$ such that
\bea
\label{inteps_eta1}
\int_{\Gamma^{k,\ell}}(p_{k\ell h}^1+\epsilon_{k,\ell}+\alpha u_{kh}^1)
\psi_{k,\ell}
= \int_{\Gamma^{k,\ell}}(-p_{\ell kh}^1+\alpha \eta_{\ell,k} +\alpha u_{\ell
h}^1) \psi_{k,\ell}
,\ \forall \psi_{k,\ell} \in \tilde W_h^{k,\ell}\hspace{10mm}
\eea
\bea
\label{inteps_eta2}
\int_{\Gamma^{k,\ell}}(p_{\ell kh}^1+\alpha \eta_{\ell,k} +\alpha u_{\ell
h}^1) \psi_{\ell,k} =
\int_{\Gamma^{k,\ell}}(-p_{k\ell h}^1-\epsilon_{k,\ell}+\alpha u_{kh}^1)
\psi_{\ell,k}
,\ \forall \psi_{\ell,k} \in \tilde W_h^{\ell,k}.\hspace{10mm}
\eea
In order to correct that, without polluting
(\ref{bestfit_u.2})-(\ref{bestfit_p.2}), for each couple $(k,\ell)$ we choose
one side, say the smaller indexed one, hereafter we shall also assume
that each couple $(k,\ell)$ is ordered by $k < \ell$. Associated to that
choice, we define $\epsilon_{k,\ell} \in \tilde W_h^{k,\ell}$,
$\eta_{\ell,k} \in {\cal Y}_h^{\ell,k} \cap H_0^1(\Gamma^{k,\ell})$,
such that $(\tilde{\uu}_h,\tilde{\up}_h)$ satisfy (\ref{disc.const})
where we define
\bea
\tilde{u}_{\ell h}=u^1_{\ell h}+\sum_{k<\ell} {\cal R}_{\ell,k}(\eta_{\ell,k})
\nonumber
\eea
\bea\label{defptilde}
\tilde{p}_{k \ell h}=p^1_{k \ell h}+ \epsilon_{k,\ell} \quad (\mbox{for } k
< \ell)
\eea
where ${\cal R}_{\ell,k}$ is a discrete lifting operator
(see \cite{Widlund}, \cite{BG}) that to any element
of ${\cal Y}_h^{\ell,k} \cap H_0^1(\Gamma^{k,\ell})$ associates a finite
element function over $\Omega^{\ell}$ that vanishes over
$\partial\Omega^{\ell}\setminus\Gamma^{k,\ell}$ and satisfies
\bea
\forall w \in {\cal Y}_h^{\ell,k} \cap H_0^1(\Gamma^{k,\ell}),
  \ ({\cal R}_{\ell,k}(w))_{| \Gamma_{k,\ell}}=w
\nonumber\\
\label{lifting}
\| {\cal R}_{\ell,k}(w) \|_{H^1(\Omega^{\ell})}
\le c \| w \|_{H^{1 \over 2}_{00}(\Gamma^{k,\ell})}
\eea
where $c$ is $h$-independent.\\
The set of equations (\ref{inteps_eta1})-(\ref{inteps_eta2}) for $\epsilon_{k,\ell}$ and $\eta_{\ell,k}$ results
in a square system of linear algebraic  equations that can be written
as follows

\bea\label{disc.const_2}
\int_{\Gamma^{k,\ell}}(\epsilon_{k,\ell}-\alpha \eta_{\ell,k})\psi_{k,\ell}
= \int_{\Gamma^{k,\ell}} e_1 \psi_{k,\ell}
,\ \forall \psi_{k,\ell} \in \tilde W_h^{k,\ell}
\eea
\bea\label{disc.const_3}
\int_{\Gamma^{k,\ell}}(\epsilon_{k,\ell}+\alpha \eta_{\ell,k})\psi_{\ell,k}
= \int_{\Gamma^{k,\ell}} e_2 \psi_{\ell,k}
,\ \forall \psi_{\ell,k} \in \tilde W_h^{\ell,k}
\eea
with
\bea\label{e1}
e_1=-p_{k\ell h}^1-p_{\ell kh}^1+\alpha(u_{\ell h}^1-u_{kh}^1)
\eea
and
\bea\label{e2}
e_2=-p_{k\ell h}^1-p_{\ell kh}^1+\alpha(u_{kh}^1-u_{\ell h}^1).
\eea
\begin{proposition}
The linear system (\ref{disc.const_2})-(\ref{disc.const_3}) is well posed.
\end{proposition}
\\
\\
{\bf Proof:}
With these notations, (\ref{disc.const_2}) yields
\bea\label{epskl}
\epsilon_{k,\ell}=
\pi_{k,\ell}(\alpha \eta_{\ell,k} +e_1)
\eea
and (\ref{disc.const_3}) yields
\bea\label{etalk}
\alpha \eta_{\ell,k}=
\pi_{\ell,k}(-\epsilon_{k,\ell} +e_2).
\eea
As (\ref{disc.const_2})-(\ref{disc.const_3}) is a square linear 
system, it suffices to
prove uniqueness for $e_1$ and $e_2$ null. From 
(\ref{epskl})-(\ref{etalk}), we get
\[
0=\eta_{\ell,k}+ \pi_{\ell,k}\pi_{k,\ell}(\eta_{\ell,k}).
\]
so that for all $\psi_{\ell,k}$ in $\tilde W_h^{k,\ell}$,
\[
0=\int_{\Gamma^{k,\ell}}(\eta_{\ell,k}+\pi_{k,\ell}(\eta_{\ell,k}))\psi_{\ell,k}.
\]
By Lemma~\ref{lem_1}, this proves that $\eta_{\ell,k}$ is zero, thus by 
(\ref{epskl}),
$\epsilon_{k,\ell}$ is zero.\\

Let us resume the proof of theorem \ref{best-fit}: By (\ref{epskl}) 
and (\ref{etalk})
we have
\bea\label{systeme.eta}
\int_{\Gamma^{k,\ell}}(\eta_{\ell,k}+\pi_{k,\ell}(\eta_{\ell,k}))\psi_{\ell,k}
= {1 \over \alpha} \int_{\Gamma^{k,\ell}} (e_2 -
\pi_{k,\ell}(e_1))\psi_{\ell,k}
,\ \forall \psi_{\ell,k} \in \tilde W_h^{\ell,k}.
\eea
%
In order to estimate $\| \tilde{p}_{k \ell h} - p_{k,\ell} \|_{H^{-{1 \over
2}}(\Gamma^{k,\ell})}$ and $\| \tilde{u}_{\ell h} -
u_{\ell}\|_{H^1(\Omega^{\ell})}$, we first estimate $\|\eta_{\ell,k}
\|_{L^2(\Gamma^{k,\ell})}$ :\\\\
from (\ref{injectif}) and (\ref{systeme.eta}) we get
\bea\label{estim.eta}
c_1\|\eta_{\ell,k} \|_{L^2(\Gamma^{k,\ell})}^2
\le {1 \over \alpha} \|e_2 - \pi_{k,\ell}(e_1)\|_{L^2(\Gamma^{k,\ell})}
\|\psi_{\ell,k}\|_{L^2(\Gamma^{k,\ell})}
\eea
and using (\ref{stable}) in (\ref{estim.eta})
\bea
\|\eta_{\ell,k} \|_{L^2(\Gamma^{k,\ell})}
\le {c_2 \over \alpha c_1} \|e_2 - \pi_{k,\ell}(e_1)\|_{L^2(\Gamma^{k,\ell})}
\nonumber
\eea
hence
\bea\label{estim.eta2}
\|\eta_{\ell,k} \|_{L^2(\Gamma^{k,\ell})}
\le {c_2 \over \alpha c_1} (\| e_2 \|_{L^2(\Gamma^{k,\ell})}
+\| e_1 \|_{L^2(\Gamma^{k,\ell})})
\eea
Now, from (\ref{e1}) and (\ref{e2}), for $i=1,2$
\bea
\|e_i \|_{L^2(\Gamma^{k,\ell})}
\le \|p_{k\ell h}^1+p_{\ell kh}^1\|_{L^2(\Gamma^{k,\ell})}
+\alpha \|u_{\ell h}^1-u_{kh}^1 \|_{L^2(\Gamma^{k,\ell})}
\nonumber
\eea
and recalling that $p_{k,\ell}=\frac{\partial u}{\partial {\bf n}_k}=
-\frac{\partial u}{\partial {\bf n}_{\ell}}=-p_{\ell,k}$
over each $\Gamma^{k,\ell}$
\bea
\|p_{k\ell h}^1+p_{\ell kh}^1\|_{L^2(\Gamma^{k,\ell})}
&\le& \|p_{k\ell h}^1-p_{k,\ell}\|_{L^2(\Gamma^{k,\ell})}
+\|p_{\ell k h}^1-p_{\ell,k}\|_{L^2(\Gamma^{k,\ell})}
\nonumber\\
\|u_{\ell h}^1-u_{kh}^1 \|_{L^2(\Gamma^{k,\ell})}
&\le& \|u_{k h}^1-u_k \|_{L^2(\Gamma^{k,\ell})}
+\|u_{\ell h}^1-u_{\ell} \|_{L^2(\Gamma^{k,\ell})}
\nonumber
\eea
so that, using (\ref{bestfit_u.2}) and (\ref{bestfit_p}), we derive
for $i=1,2$ and $0\le m\le M-1$
\bea\label{estim.ei}
\| e_i \|_{L^2(\Gamma^{k,\ell})}
\le c \alpha h^{{3 \over 2}+m}
(\|u_k\|_{H^{2+m}(\Omega^k)}+\|u_{\ell}\|_{H^{2+m}(\Omega^{\ell})})
+ch^{{1 \over 2}+m}
\| p_{k,\ell} \|_{H^{{1 \over 2}+m}(\Gamma^{k,\ell})}\hspace{10mm}
\eea
and (\ref{estim.eta2}) yields for $0\le m\le M-1$
\bea\label{estim.eta3}
\|\eta_{\ell,k} \|_{L^2(\Gamma^{k,\ell})}
\le c h^{{3 \over 2}+m}
(\|u_k\|_{H^{2+m}(\Omega^k)}+\|u_{\ell}\|_{H^{2+m}(\Omega^{\ell})})
+{ch^{{1 \over 2}+m} \over \alpha}
\| p_{k,\ell} \|_{H^{{1 \over 2}+m}(\Gamma^{k,\ell})}.\hspace{10mm}
\eea
%
We can now evaluate $\| \tilde{p}_{k \ell h} - p_{k,\ell} \|_{H^{-{1 \over
2}}(\Gamma^{k,\ell})}$, using (\ref{defptilde}) :
\bea\label{eval1}
\| \tilde{p}_{k \ell h} - p_{k,\ell} \|_{H^{-{1 \over 2}}(\Gamma^{k,\ell})}
\le \| \epsilon_{k,\ell} \|_{H^{-{1 \over 2}}(\Gamma^{k,\ell})}
+ \| p_{k \ell h}^1 - p_{k,\ell} \|_{H^{-{1 \over 2}}(\Gamma^{k,\ell})}.
\eea
%
The term $\| p_{k \ell h}^1 - p_{k,\ell} \|_{H^{-{1 \over
2}}(\Gamma^{k,\ell})}$ is estimated in (\ref{bestfit_p.2}), so let us focus
on the term
$\| \epsilon_{k,\ell} \|_{H^{-{1 \over 2}}(\Gamma^{k,\ell})}$.
 From (\ref{epskl}) we have,
\bea\label{majoeps}
  \| \epsilon_{k,\ell} \|_{H^{-{1 \over 2}}(\Gamma^{k,\ell})}
\le \alpha\|\eta_{\ell,k} \|_{H^{-{1 \over 2}}(\Gamma^{k,\ell})}
+\|e_1 \|_{H^{-{1 \over 2}}(\Gamma^{k,\ell})}
+ \|(Id-\pi_{k,\ell})(\alpha \eta_{\ell,k} +e_1) \|_{H^{-{1 \over
2}}(\Gamma^{k,\ell})}.\hspace{10mm}
\eea
To evaluate $\|e_1 \|_{H^{-{1 \over 2}}(\Gamma^{k,\ell})}$ we proceed
as for $\|e_1 \|_{L^2(\Gamma^{k,\ell})}$ and from (\ref{bestfit_u.3}) and
(\ref{bestfit_p.2}) we have, for $i=1,2$ and $0\le m\le M$:
\bea\label{estim.ei.2}
\|e_i \|_{H^{-{1 \over 2}}(\Gamma^{k,\ell})}
\le c \alpha h^{2+m}
(\|u_k\|_{H^{2+m}(\Omega^k)}+\|u_{\ell}\|_{H^{2+m}(\Omega^{\ell})})
+ch^{1+m} \| p_{k,\ell} \|_{H^{{1 \over 2}+m}(\Gamma^{k,\ell})}.
\hspace{10mm}
\eea
The third term in the right member of (\ref{majoeps}) satisfies
\bea
\|(Id-\pi_{k,\ell})(\alpha \eta_{\ell,k} +e_1) \|_{H^{-{1 \over
2}}(\Gamma^{k,\ell})} \le c \sqrt{h} \|\alpha \eta_{\ell,k} +e_1
\|_{L^2(\Gamma^{k,\ell})}. \nonumber
\eea
Then, using (\ref{estim.eta3}) and (\ref{estim.ei}) yields
\bea
\|(Id-\pi_{k,\ell})(\alpha \eta_{\ell,k} +e_1) \|_{H^{-{1 \over
2}}(\Gamma^{k,\ell})}
&\le& c\alpha h^{2+m}
(\|u_k\|_{H^{2+m}(\Omega^k)}+\|u_{\ell}\|_{H^{2+m}(\Omega^{\ell})})\nonumber\\&
&+ \ c h^{1+m} \| p_{k,\ell} \|_{H^{{1 \over 2}+m}(\Gamma^{k,\ell})}.
\nonumber
\eea
In order to estimate the term $\|\eta_{\ell,k} \|_{H^{-{1 \over
2}}(\Gamma^{k,\ell})}$ in (\ref{majoeps}), we
use (\ref{systeme.eta}):
\bea
2\int_{\Gamma^{k,\ell}} \eta_{\ell,k} \psi_{\ell,k}&=&
2\int_{\Gamma^{k,\ell}} \eta_{\ell,k} \psi_{\ell,k}\nonumber\\
&-&(
\int_{\Gamma^{k,\ell}}(\eta_{\ell,k}+\pi_{k,\ell}(\eta_{\ell,k}))
\psi_{\ell,k}
-{1 \over \alpha} \int_{\Gamma^{k,\ell}} (e_2 -
\pi_{k,\ell}(e_1))
\psi_{\ell,k}
)\nonumber
\eea
that is
\bea
2\int_{\Gamma^{k,\ell}} \eta_{\ell,k} \psi_{\ell,k} = \int_{\Gamma^{k,\ell}}
(\eta_{\ell,k}-\pi_{k,\ell}\eta_{\ell,k})\psi_{\ell,k}
+{1 \over \alpha} \int_{\Gamma^{k,\ell}} (e_2 -
\pi_{k,\ell}(e_1)) \psi_{\ell,k}.
\nonumber
\eea
Using the symmetry of the operator $\pi_{k,\ell}$ we deduce
\bea
2\int_{\Gamma^{k,\ell}} \eta_{\ell,k} \psi_{\ell,k} = \int_{\Gamma^{k,\ell}}
(\psi_{\ell,k}-\pi_{k,\ell}\psi_{\ell,k})   \eta_{\ell,k}
+{1 \over \alpha} \int_{\Gamma^{k,\ell}} (e_2 -
\pi_{k,\ell}(e_1)) \psi_{\ell,k}.
\nonumber
\eea
then, from (\ref{eq:propr-pilk}) yields
\bea
|\int_{\Gamma^{k,\ell}} \eta_{\ell,k} \psi_{\ell,k}|
\le c\sqrt{h}\|\eta_{\ell,k}\|_{L^2(\Gamma^{k,\ell})}
\|\psi_{\ell,k}\|_{H^{1 \over 2}(\Gamma^{k,\ell})}
+{1\over \alpha}\|e_2 - \pi_{k,\ell}(e_1)\|_{H^{-{1 \over 2}}(\Gamma^{k,\ell})}
\|\psi_{\ell,k}\|_{H^{1 \over 2}(\Gamma^{k,\ell})}
\nonumber
\eea
and thus,
  we have
\bea
\|\eta_{\ell,k}\|_{H^{-{1 \over 2}}(\Gamma^{k,\ell})}
\le c \sqrt{h}\|\eta_{\ell,k}\|_{L^2(\Gamma^{k,\ell})}+ {c\over \alpha}
\|e_2 - \pi_{k,\ell}(e_1)\|_{H^{-{1
\over 2}}(\Gamma^{k,\ell})}. \nonumber
\eea
Then, using (\ref{estim.eta3}) and the fact that
\bea
\|e_2 - \pi_{k,\ell}(e_1)\|_{H^{-{1 \over 2}}(\Gamma^{k,\ell})}
\le \|e_2\|_{H^{-{1 \over 2}}(\Gamma^{k,\ell})}
+\|e_1\|_{H^{-{1 \over 2}}(\Gamma^{k,\ell})}
+\|e_1- \pi_{k,\ell}(e_1)\|_{H^{-{1 \over 2}}(\Gamma^{k,\ell})} \nonumber\\
\le \|e_2\|_{H^{-{1 \over 2}}(\Gamma^{k,\ell})}
+\|e_1\|_{H^{-{1 \over 2}}(\Gamma^{k,\ell})}
+c\sqrt{h} \|e_1\|_{L^2(\Gamma^{k,\ell})} \nonumber
\eea
with (\ref{estim.ei}) and (\ref{estim.ei.2}) yields for $0\le m \le M-1$
\bea
\|\eta_{\ell,k}\|_{H^{-{1 \over 2}}(\Gamma^{k,\ell})}
\le c h^{2+m}
(\|u_k\|_{H^{2+m}(\Omega^k)}+\|u_{\ell}\|_{H^{2+m}(\Omega^{\ell})})
+{ch^{1+m} \over \alpha}
\| p_{k,\ell} \|_{H^{{1 \over 2}+m}(\Gamma^{k,\ell})}. \nonumber
\eea
Using the previous inequality in (\ref{majoeps}), (\ref{eval1}) yields
\bea\label{estim.p}
\| \tilde{p}_{k \ell h} - p_{k,\ell} \|_{H^{-{1 \over 2}}(\Gamma^{k,\ell})}
\le c\alpha h^{2+m}
(\|u_k\|_{H^{2+m}(\Omega^k)}+\|u_{\ell}\|_{H^{2+m}(\Omega^{\ell})})
+c h^{1+m} \| p_{k,\ell} \|_{H^{{1 \over 2}+m}(\Gamma^{k,\ell})}.\hspace{10mm}
\eea
Let us now estimate $\| \tilde{u}_{\ell h} - u_{\ell}\|_{H^1(\Omega^{\ell})}$ :
\bea\label{estim.u2}
\| \tilde{u}_{\ell h} - u_{\ell}\|_{H^1(\Omega^{\ell})}
&\le& \| u_{\ell h}^1 - u_{\ell}\|_{H^1(\Omega^{\ell})}
+ \sum_{k < \ell} \| {\cal R}_{\ell,k}(\eta_{\ell,k}) \|_{H^1(\Omega^{\ell})}
\eea
and from (\ref{lifting})
\bea
\| {\cal R}_{\ell,k}(\eta_{\ell,k}) \|_{H^1(\Omega^{\ell})}
\le c \|\eta_{\ell,k} \|_{H^{1 \over 2}_{00}(\Gamma^{k,\ell})} \nonumber
\eea
then, with an inverse inequality
\bea
\| {\cal R}_{\ell,k}(\eta_{\ell,k}) \|_{H^1(\Omega^{\ell})}
\le c h^{- {1 \over 2}} \|\eta_{\ell,k} \|_{L^2(\Gamma^{k,\ell})}. \nonumber
\eea
Hence, from (\ref{estim.eta3}) we have for $0\le m \le M-1$
\bea
\| {\cal R}_{\ell,k}(\eta_{\ell,k}) \|_{H^1(\Omega^{\ell})}
\le
ch^{1+m}(\|u_k\|_{H^{2+m}(\Omega^k)}+\|u_{\ell}\|_{H^{2+m}(\Omega^{\ell})})
\nonumber\\
+{c h^m\over \alpha}  \| p_{k,\ell} \|_{H^{{1 \over 2}+m}(\Gamma^{k,\ell})}
\nonumber
\eea
%
and (\ref{estim.u2}) yields
\bea\label{estim.u3}
\| \tilde{u}_{\ell h} - u_{\ell}\|_{H^1(\Omega^{\ell})}
\le ch^{1+m} \|u_{\ell}\|_{H^{2+m}(\Omega^{\ell})}
+ ch^{1+m}\sum_{k < \ell} \|u_k\|_{H^{2+m}(\Omega^k)}
\nonumber\\
+{c h^m \over \alpha} \sum_{k < \ell} \| p_{k,\ell} \|_{H^{{1 \over
2}+m}(\Gamma^{k,\ell})}.
\eea
As $u \in H^{2+m}(\Omega^{\ell})$,
\bea
\| p_{\ell} \|_{H^{{1 \over 2}+m}(\partial \Omega^{\ell})}
\le c \| u_{\ell} \|_{H^{2+m}(\Omega^{\ell})} \nonumber
\eea
and for $0\le m \le M-1$
\bea
\| \tilde{\uu}_h -\uu\|_*
\le {c h^m\over \alpha} \sum_{i=1}^K \| \uu\|_{H^{2+m}(\Omega^k)}. \nonumber
\eea
%
\\\\
{\bf Proof of theorem \ref{best-fit.2}:} The proof is the same that
for theorem \ref{best-fit}, except that the relation (\ref{bestfit_p})
for $0\le m \le M-1$
is changed using the following lemma

\begin{lem}\label{lem.logh}
%
$p_{k \ell h}^1$ satisfy for $0\le m \le M-1$ the error estimate
\bea\label{eq:logh}
\| p_{k \ell h}^1 - p_{k,\ell} \|_{L^2(\Gamma^{k,\ell})}
&\le& c h^{{3 \over 2}+m}\ (\log h)^{\beta(m)} \| p_{k,\ell} \|_{H^{{3 \over
2}+m}(\Gamma^{k,\ell})}.
\eea
\end{lem}
%
Therefore, (\ref{bestfit_p.2}) is changed in
\bea
\| p_{k \ell h}^1 - p_{k,\ell} \|_{H^{-{1 \over 2}}(\Gamma^{k,\ell})}
&\le& c h^{2+m}\ (\log h)^{\beta(m)}\| p_{k,\ell} \|_{H^{{3 \over
2}+m}(\Gamma^{k,\ell})}
\nonumber
\eea
and (\ref{estim.p}) is changed in
\bea
\| \tilde{p}_{k \ell h} - p_{k,\ell} \|_{H^{-{1 \over 2}}(\Gamma^{k,\ell})}
\le c\alpha h^{2+m}
(\|u_k\|_{H^{2+m}(\Omega^k)}+\|u_{\ell}\|_{H^{2+m}(\Omega^{\ell})})\nonumber \\
  +c h^{2+m} \ (\log h)^{\beta(m)} \| p_{k,\ell} \|_{H^{{3\over
2}+m}(\Gamma^{k,\ell})} \nonumber
\eea
and (\ref{estim.u3}) is changed in
\bea
\| \tilde{u}_{\ell h} - u_{\ell}\|_{H^1(\Omega^{\ell})}
\le c h^{1+m} \|u_{\ell}\|_{H^{2+m}(\Omega^{\ell})}
+ c h^{1+m} \sum_{k < \ell} \|u_k\|_{H^{2+m}(\Omega^k)}\nonumber \\
+{c h^{1+m} \over \alpha} \ (\log h)^{\beta(m)} \sum_{k < \ell} \| p_{k,\ell}
\|_{H^{{3 \over 2}+m}(\Gamma^{k,\ell})}. \nonumber
\eea
%
{\bf Proof of lemma \ref{lem.logh}:}
For $0\le m < M-1, \ \beta(m)=0$ and the estimate (\ref{eq:logh}) is standard.
For $m=M-1$, let ${\bar p}_{k \ell h}$ be the unique element of $\tilde W_h^{k,\ell}$ defined as follows :
\begin{itemize}
\item
$({\bar p}_{k \ell h})_{|[x_1^{\ell,k},x_{n-1}^{\ell,k}]}$ coincide
with the interpolate of degree $M$ of $p_{k,\ell}$.
\item
$({\bar p}_{k \ell h})_{|[x_0^{\ell,k},x_{1}^{\ell,k}]}$ and 
$({\bar p}_{k \ell h})_{|[x_{n-1}^{\ell,k},x_{n}^{\ell,k}]}$ coincide
with the interpolate of degree $M-1$ of $p_{k,\ell}$.
\end{itemize}
Then, we have
\bea
\| p_{k \ell h}^1 - p_{k,\ell} \|^2_{L^2(\Gamma^{k,\ell})}
\le \| {\bar p}_{k \ell h} - p_{k,\ell} \|^2_{L^2(\Gamma^{k,\ell})}.\nonumber
\eea
Using Deny-Lions theorem we have
\bea
\| p_{k \ell h}^1 - p_{k,\ell} \|^2_{L^2(\Gamma^{k,\ell})}
\le h^{2M}\int_{x_0^{\ell,k}}^{x_1^{\ell,k}}|\frac{d^M p_{k,\ell}}{dx^M}|^2
+h^{2(M+{1 \over 2})}\|p_{k,\ell}\|_{H^{M+{1 \over 2}}([x_1^{\ell,k},x_{n-1}^{\ell,k}])}
+h^{2M}\int_{x_{n-1}^{\ell,k}}^{x_n^{\ell,k}}|\frac{d^M p_{k,\ell}}{dx^M}|^2.
\nonumber
\eea
Let $\varphi=\frac{d^M p_{k,\ell}}{dx^M}$.
In order to analyse the two extreme contributions, we use
H\"older's inequality : 
\bea
\int_0^h\varphi^2 \le h^{1-{1 \over p}}\|\varphi\|^2_{L^p(0,h)}.\nonumber
\eea
Then, we use the estimate
\bea
\|\varphi\|_{L^p(0,h)} \le cp\|\varphi\|_{H^{1 \over 2}(0,h)},\nonumber
\eea
where $c$ is a constant. Thus we have
\bea
\int_0^h\varphi^2 \le cp^2 h^{1-{1 \over p}}\|\varphi\|^2_{H^{1 \over 2}(0,h)}.\nonumber
\eea
Then, 
\bea
h^{2M}\int_{x_0^{\ell,k}}^{x_1^{\ell,k}}|\frac{d^M p_{k,\ell}}{dx^M}|^2 \le c h^{2M}p^2 h^{1-{1 \over p}}\|\frac{d^M p_{k,\ell}}{dx^M}\|^2_{H^{1 \over 2}(x_0^{\ell,k},x_1^{\ell,k})}.
\nonumber
\eea
Now we take $p=\log h$ and thus we obtain
\bea
h^{2M}\int_{x_0^{\ell,k}}^{x_1^{\ell,k}}|\frac{d^M p_{k,\ell}}{dx^M}|^2 \le c (h^{M+{1 \over 2}} \log h)^2\|p_{k,\ell}\|^2_{H^{{1 \over 2}+M}(x_0^{\ell,k},x_1^{\ell,k})}.
\nonumber
\eea
In a same way we have
\bea
h^{2M}\int_{x_{n-1}^{\ell,k}}^{x_n^{\ell,k}}|\frac{d^M p_{k,\ell}}{dx^M}|^2 \le c (h^{M+{1 \over 2}} \log h)^2\|p_{k,\ell}\|^2_{H^{{1 \over 2}+M}(x_{n-1}^{\ell,k},x_n^{\ell,k})},
\nonumber
\eea
and thus we obtain
\bea
\| p_{k \ell h}^1 - p_{k,\ell} \|_{L^2(\Gamma^{k,\ell})}
\le  c (h^{M+{1 \over 2}} \log h)\|p_{k,\ell}\|_{H^{{1 \over 2}+M}(\Gamma^{k,\ell})},\nonumber
\eea
which ends the proof of lemma \ref{lem.logh}.
\subsection{Error Estimates.}
Thanks to (\ref{estimuuh}), we have the following error estimates:
\begin{theo}
\label{error-estimate}
Assume that the solution $u$ of
(\ref{initial_BVP1})-(\ref{initial_BVP2}) is in $H^2(\Omega)\cap H^1_0(\Omega)$, and $u_k=u_{|\Omega^k}\in
H^{2+m}(\Omega^k)$, with
$M-1\ge m \ge 0$,
and let
$p_{k,\ell}=\frac{\partial u}{\partial {\bf n}_k}$
over each $\Gamma^{k,\ell}$.
Then, there exists a constant $c$ independent of $h$ and $\alpha$
such that
\bea
\| \uu_h -\uu\|_* + \| \up_h - \up \|_{-{1 \over 2},*}
\le c(\alpha h^{2+m}+h^{1+m} )
\sum_{k=1}^K \| \uu\|_{H^{2+m}(\Omega^k)} \nonumber\\
+ \ c ({h^m\over \alpha} + h^{1+m})\sum_{k=1}^K \sum_{\ell}\| 
p_{k,\ell} \|_{H^{{1
\over 2}+m}(\Gamma^{k,\ell})}.
\nonumber
\eea
%
\end{theo}
%
\begin{theo}
\label{error-estimate2}
Assume that the solution $u$ of
(\ref{initial_BVP1})-(\ref{initial_BVP2}) is in $H^2(\Omega)\cap H^1_0(\Omega)$, $u_k=u_{|\Omega^k}\in
H^{2+m}(\Omega^k)$,
and
$p_{k,\ell}=\frac{\partial u}{\partial {\bf n}_k}$
is in $H^{{3 \over 2}+m}(\Gamma_{k,\ell})$ with $M-1\ge m \ge 0$.
Then there exists a constant $c$ independent of $h$ and $\alpha$
such that
\bea
\| \uu_h -\uu\|_* + \| \up_h - \up\|_{-{1 \over 2},*}
\le c(\alpha h^{2+m}+h^{1+m} )
\sum_{k=1}^K \| \uu\|_{H^{2+m}(\Omega^k)} \nonumber\\
+ \ c ({h^{1+m}\over \alpha} + h^{2+m}) (\log h)^{\beta(m)} 
\sum_{k=1}^K \sum_{\ell}\|
p_{k,\ell}
\|_{H^{{3
\over 2}+m}(\Gamma^{k,\ell})}
\nonumber
\eea
with $\beta(m)=0$ if
$m\le M-2$ and $\beta(m)=1$ if $m=M-1$.
\end{theo}
%
%
\begin{rem}
Let us consider a $P_1$ finite element approximation. If the solution $\uu$
of (\ref{initial_BVP1})-(\ref{initial_BVP2}) is in
$\prod_{k=1}^K H^2_*(\Omega^k)$ and
$\alpha$ is a constant independent of $h$ then, from theorem
\ref{error-estimate},
\bea
\| \tilde{\uu}_h -\uu\|_*
\le c \sum_{i=1}^K \| \uu\|_{H^2(\Omega^k)} \nonumber
\eea
and this result is not optimal. In order to improve this fact, we have to
choose a parameter $\alpha$ which depends on $h$, or assume that
$\uu \in \prod_{k=1}^K H^3_*(\Omega^k)$, or assume that
$\uu \in \prod_{k=1}^K H^2_*(\Omega^k)$ and $p_{k,\ell}=\frac{\partial
u}{\partial {\bf n}_k} \in H^{3 \over 2}(\Gamma_{k,\ell})$:
\begin{itemize}
\item
If  the solution $\uu$ of (\ref{initial_BVP1})-(\ref{initial_BVP2}) is in
$\prod_{k=1}^K H^2_*(\Omega^k)$ and
$\alpha={c \over h}$, then
\bea
\| \tilde{\uu}_h -\uu\|_*
\le c h\sum_{i=1}^K \| \uu\|_{H^2(\Omega^k)} \nonumber
\eea
\item
If  the solution $\uu$ of (\ref{initial_BVP1})-(\ref{initial_BVP2}) is in
$\prod_{k=1}^K H^2_*(\Omega^k)$, $p_{k,\ell}=\frac{\partial u}{\partial
{\bf n}_k}$ is in $H^{3 \over 2}(\Gamma_{k,\ell})$
  and $\alpha$ is a constant independent of
$h$ then
\bea
\| \tilde{\uu}_h -\uu\|_* = O(h |\log(h)|)
. \nonumber
\eea
\end{itemize}
\end{rem}
%
\subsection{Analysis of the best fit in 3D}
\begin{figure}[H]
\centering 
  \psfrag{boun}{Boundary of $\displaystyle \Gamma $} 
  \psfrag{T1}{$\displaystyle T_{i-1,r}$}
  \psfrag{T2}{$\displaystyle T_{i,c}$}
  \psfrag{T3}{$\displaystyle T_{i,r}$}
  \psfrag{T4}{$\displaystyle T_{j',c}$} 
  \psfrag{T5}{$\displaystyle T_{c'}$}
  \psfrag{a}{$\displaystyle a$}
  \psfrag{b}{$\displaystyle b$}
  \psfrag{c}{$\displaystyle c$}
  \psfrag{cp}{$\displaystyle c'$}
  \psfrag{vi}{$\displaystyle v_i$}
  \psfrag{vip}{$\displaystyle v_{i+1}$}
  \psfrag{vim}{$\displaystyle v_{i-1}$}
  \psfrag{ui}{$\displaystyle u_i$}
  \psfrag{uip}{$\displaystyle u_{i+1}$}
  \includegraphics[width=9cm]{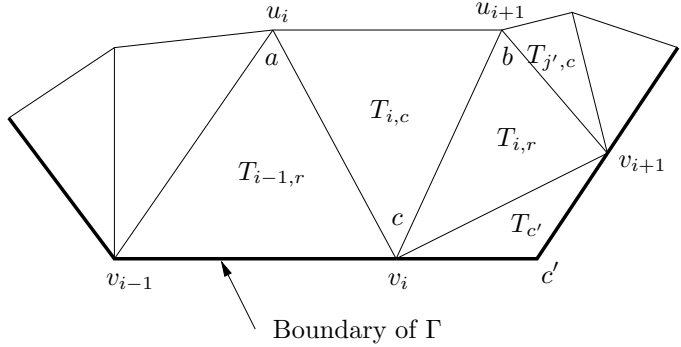}
  \caption{Two different situations of 2D triangulation of the interface
  $\Gamma$, next to it's boundary (near cross points)}
  \label{fig:triangles}
  \end{figure}

 In this section, we prove lemma \ref{lem_1} for a $P_1$-discretization in 3D. We shall use the construction proposed in \cite{BraessDahmen}. 
 In order to make the reading easy, we shall recall the notations of the above mentioned paper. The analysis is done on one subdomain $\Omega^k$ that will
 be fixed in what follows. A typical interface between this subdomain and a generic subdomain $\Omega_l$ will be denoted by $\Gamma$. We denote by $\T$ the restriction to
 $\Gamma$ 
 of the triangulation ${\cal T}_h^k$. Let $S({\T})$ denote the space of piecewise linear functions with respect to $\T$ which are continuous on $\Gamma$ and vanish on its boundary. The space of the Lagrange multipliers on $\Gamma$, defined below, will be denoted by $M({\mathcal T})$.
 In 2D, the requirement dim $M({\mathcal T})=$dim $S({\T})$ can be satisfied
 by lowering the degree of the finite elements on the intervals next
 to the end points of the interface. In 3D, it is slightly more complex (see
 \cite{BBM}). Thus, we shall use the construction proposed in \cite{BraessDahmen} in the case where all the vertices of the boundary of $\Gamma$ are connected to zero or two vertices in the interior of $\Gamma$ (figure~\ref{fig:triangles}). Let ${\mathcal V}$, ${\mathcal V}_0$, 
 $\partial {\mathcal V}$ denote respectively the set of all the vertices
 of ${\mathcal T}$, the vertices in the interior of $\Gamma$, and the vertices
 on the boundary of $\Gamma$. The finite element basis functions will
 be denoted by $\Phi_a, \ a \in {\mathcal V}$. Thus,
 \bea
 S({\T})=\mbox{ span }\{\Phi_a : a \in {\mathcal V}_0\}. \nonumber
 \eea
 For $a \in {\mathcal V}$, let $\sigma_a$ denote the support of $\Phi_a$,
 \bea
 \sigma_a:=\bigcup\{T \in \T :  a \in T\}, \nonumber
 \eea
 and let ${\mathcal N}_a$ be the set of neighboring vertices in ${\mathcal V}_0$ of a:
 \bea
 {\mathcal N}_a:=\{b \in {\mathcal V}_0 :  b \in \sigma_a\}.
 \nonumber
 \eea
 Thus,
 \bea
 {\mathcal N}=\ds \bigcup_{ a \in \partial{\mathcal V}} {\mathcal N}_a \nonumber
 \eea
 is the set of those interior vertices which have a neighbor on the boundary
 of $\Gamma$. If some triangle $T \in \T$ has all its vertices
 on the boundary of $\Gamma$, then there exists one (corner) vertex
 which has no neighbor in ${\mathcal V}_0 $. Let $\T_c$ be the set of
 triangles $T \in \T$ which have all their vertices
 on the boundary of $\Gamma$. For $T \in \T_c$, we denote by $c_T$ the only
 vertex of $T$ that has no interior neighbour 
 (such a vertex is unique as soon as the triangulation is fine enough).
 Let ${\mathcal N}_c$ denote the vertices $a_T$ of ${\mathcal N}$ which belong to a triangle adjacent to a triangle $T \in \T_c$.
 Now, we define
 the space $M({\mathcal T})$ by
 \bea
 M({\mathcal T}):=\mbox{ span }\{\hat{\Phi}_a, \ a \in {\mathcal V}_0 \}, \nonumber
 \eea
 where the basis functions $\hat{\Phi}_a$ are defined as follows :
 \bea
 \hat{\Phi}_a:=\left\{\begin{array}{ll}
 \Phi_a, & a \in {\mathcal V}_0 \setminus {\mathcal N}\\
 \Phi_a+\ds\sum_{b \in \partial{\mathcal V} \cap \sigma_a} A_{b,a}\Phi_b &
 a \in {\mathcal N} \setminus {\mathcal N}_c\\
 \Phi_{a_T}+\ds\sum_{b \in \partial{\mathcal V} \cap \sigma_{a_T}} A_{b,{a_T}}\Phi_b 
 +\Phi_{c_T} & a={a_T} \in {\mathcal N}_c
 \end{array}\right.\nonumber
 \eea
the weights $A_{b,a}$ being defined in (\ref{eq:FirstABbraess}). For all boundary nodes $c \in \partial{\mathcal V}$ connected to two interior nodes $a$ and $b$, if $T_a$ (resp. $T_b$) denote the triangle having an edge on $\partial\Gamma$ and $a$ (resp. $b$) as the opposite vertex, then
the weights are defined such that (see \cite{BraessDahmen})
\begin{equation}\label{eq:FirstABbraess}
A_{c,a}+A_{c,b}=1 \text{ and } |T_b|A_{c,a}=|T_a|A_{c,b}.
\end{equation}
$M({\mathcal T})$ is the notation introduced in \cite{BraessDahmen}, that
we use here for the sake of clarity. Corresponding to our previous
notation, $M({\mathcal T}) \equiv \tilde W_h^{k,\ell}$. \\To any $u\in S({\mathcal T})$, $u=\sum_{a\in {\mathcal V}_0} u(a)  \Phi_a$, we associate $v\in M({\mathcal T})$ where $v=\sum_{a\in {\mathcal V}_0} u(a)   \hat{\Phi}_a$. More explicitly, that means that
to any $u \in S({\T})$, we associate an element
$v \in M({\mathcal T})$ as follows (see figure~\ref{fig:triangles}): 
\begin{enumerate}[(i)]
\item $v$ is a piecewise linear finite element on $\T$ 
\item for all interior nodes $a$, $v(a):=u(a)$
\item for all boundary nodes $c$, by assumption we have two situations:
\subitem - $c$ is connected to two interior nodes denoted by $a$ and $b$. \\Then, $v(c):=A u(a)+ B u(b)$ where 
\begin{equation}\label{eq:ABbraess}
A+B=1 \text{ and } |T_b|A=|T_a|B
\end{equation}
 where $T_a$ (resp. $T_b$) is the triangle having an edge on $\partial\Gamma$ and $a$ (resp. $b$) as the opposite vertex. 
\subitem - $c$ is not connected to any interior point. We consider the triangle adjacent to the triangle to which $c$ belongs to. This triangle has one interior node denoted by $b$. Then, we define $v(b):=u(b)$. 
\end{enumerate}
Then, using the uniform regularity of $\T$, it is easy to check that there exists a constant $c$ independent of $h$ such that
\begin{equation}
\|v\|_{L^2(\Gamma)} \le c\|u\|_{L^2(\Gamma)}. \nonumber
\end{equation}
We shall need the following technical assumption:\\
{\bf Assumption} {\em 
Let  $0 < C \le 2/3$. For any triangle $T_{c'}$ having all three vertices on the boundary of $\T$ (see figure~\ref{fig:triangles}), we consider the two triangles $T_{i,c}$ and $T_{j',c}$ surrounding $T_{c'}$. We assume that
\[
\frac{1}{24}\min(|T_{i,c}|,|T_{j',c}|) > \frac{C}{2}|T_{c'}|.
\]
}

In order to prove lemma \ref{lem_1}, we prove the following lemma:
\begin{lem}\label{lem.3D}
Let  $0 < C \le 2/3$, we assume the above assumption and that $\T$ is uniformly regular. Let $u \in S({\T})$ and let $v \in M({\mathcal T})$ constructed from $u$ as
explained above ((i)-(iii)).\\
Then, there exists $c>0$ such that,
\begin{equation}\label{eq:stabesti}
\int_{\Gamma}(uv-\frac{C}{2}(u-v)^2) \ge c \int_{\Gamma}u^2.
\end{equation}
\end{lem}
%
{\bf Proof of lemma \ref{lem.3D}:}
Let us introduce the notation
\begin{equation}
Q_{\Gamma}:=\int_{\Gamma}(uv-\frac{C}{2}(u-v)^2).  \nonumber
\end{equation}
We have
\[
Q_{\Gamma}=\frac{1}{4}\int_{\Gamma}(u+v)^2-(1+2C)(u-v)^2.
\]
In order to estimate $Q_{\Gamma}$, we remark that 
\[
Q_{\Gamma}=\sum_{T\in\T} Q_T
\]
where 
\[
 Q_T=\frac{1}{4}\int_{T}(u+v)^2-(1+2C)(u-v)^2.
\]
We have four kinds of triangles:
\begin{enumerate}
\item Inner triangles i.e they don't touch the boundary of $\Gamma$.
\item Triangles which have only one vertex on the boundary
\item Triangles which have two vertices on the boundary
\item Triangles which have three vertices on the boundary
\end{enumerate}

\paragraph{Inner triangles}
On an inner triangle $T$, $u=v$ so that for all $C>0$, we have
\[
 Q_T \ge c\int_T u^2
\]
for  $c\le 1$.

\paragraph{Triangles having only one vertex on the boundary}
Let $T_{i,c}$ be such a triangle (see figure~\ref{fig:triangles}).We use the following notations: $u_i=u(a)=v(a)$, $u_{i+1}=u(b)=v(b)$ and $v_i=v(c)$. First notice that we have (remember $u(c)=0$)
\be\label{eq:intu2Tic}
\int_{T_{i,c}}u^2 = \frac{|T_{i,c}|}{12}\biggl(
u_{i}^2+u_{i+1}^2 +(u_{i}+u_{i+1})^2\biggl)
= \frac{|T_{i,c}|}{12}\biggl(
2u_{i}^2+2u_{i+1}^2 +2u_{i}u_{i+1}\biggl) \nonumber
\ee
see for example \cite{Braess} (II.8.4). As for $Q_{T_{i,c}}$, we have
\bea
Q_{T_{i,c}}=\frac{|T_{i,c}|}{48}\biggl
( (2u_{i})^2+(2u_{i+1})^2+(Au_{i}+Bu_{i+1})^2 \nonumber \\ 
+(2u_{i}+2u_{i+1}+Au_{i}+Bu_{i+1})^2 - 2(1+2C) (Au_{i}+Bu_{i+1})^2
\biggl ) \nonumber\\
=\frac{|T_{i,c}|}{48}\biggl
( 8u_{i}^2+8u_{i+1}^2+8u_{i}u_{i+1}+4(u_{i}+u_{i+1})(Au_{i}+Bu_{i+1}) \nonumber\\
-4C (Au_{i}+Bu_{i+1})^2.
\biggl ) \nonumber
\eea
If we take $C=1$ and use $A+B=1$, we get:
\bea
Q_{T_{i,c}}&=&\frac{|T_{i,c}|}{48}\biggl(
4u_{i}^2+4u_{i+1}^2+4u_{i}u_{i+1}+4AB(u_{i}-u_{i+1})^2+4(u_{i}+u_{i+1})^2
\biggl )\nonumber\\
&\ge& \frac{1}{2}\int_{T_{i,c}}u^2. \nonumber
\eea
Hence, for all $0<C\le 1$, we have:
\be
Q_{T_{i,c}}\ge \frac{1}{2}\int_{T_{i,c}}u^2. \nonumber
\ee
Therefore,
\be
Q_{T_{i,c}}\ge c\int_{T_{i,c}}u^2 \nonumber
\ee 
 for $0<C\le 1$ and $0<c\le 1/2$.
We shall also use in the sequel the estimate:
\begin{equation}
  \label{eq:Qtic}
  Q_{T_{i,c}} \ge \frac{|T_{i,c}|}{24} u_{i+1}^2. 
\end{equation}

\paragraph{Triangles having two vertices on the boundary}
We consider now a triangle $T_{i,r}$ having two vertices on the boundary of the face $\Gamma$,  see figure~\ref{fig:triangles}. 
Let ${\cal N}_r=\{i,\ T_{i,r}$ has two vertices on the boundary of $\Gamma\}$.
First notice that we have
\be 
\int_{T_{i,r}}u^2 = \frac{|T_{i,c}|}{12}\biggl(2u_{i+1}^2\biggl). \nonumber
\ee
And we have
\bea
Q_{T_{i,r}}&=&\frac{|T_{i,r}|}{48}\biggl( 
4u_{i+1}^2+v_{i}^2+v_{i+1}^2+(2u_{i+1}+v_{i}+v_{i+1})^2 \nonumber \\
&&\ \ \ \ \ \ \ \ \ \ -(1+2C)(v_{i}^2+v_{i+1}^2+(v_{i}+v_{i+1})^2) \biggl )\nonumber \\
&=&\frac{|T_{i,r}|}{48}\biggl(
8u_{i+1}^2-4Cv_{i}^2-4Cv_{i+1}^2-4Cv_{i}v_{i+1}+4u_{i+1}(v_{i}+v_{i+1})
\biggl ). \nonumber
\eea
Then,
\bea
Q_{T_{i,r}}\ge \frac{|T_{i,r}|}{48}\biggl(
8u_{i+1}^2-6Cv_{i}^2-6Cv_{i+1}^2+4u_{i+1}(v_{i}+v_{i+1})\biggl ). \nonumber
\eea
Defining $E_i:=u_{i+1}v_i$ and $F_i:=u_{i+1}v_{i+1}$ (cf. \cite{BraessDahmen} page~11), we have:
\bea
Q_{T_{i,r}}\ge \int_{T_{i,r}}u^2+\frac{|T_{i,r}|}{48}\biggl(
-6Cv_{i}^2+4E_{i}-6Cv_{i+1}^2+4F_{i} 
 \biggl ).  \nonumber
\eea
Now we sum these terms over all the triangles having two vertices on the boundary $\Gamma$.
\begin{equation}
  \label{eq:sumqtir}
  \begin{array}{l}
\ds  \sum_{i \in {\cal N}_r} Q_{T_{i,r}} \ge \int_{\cup_{i \in {\cal N}_r} T_{i,r}} u^2 + \sum_{i \in {\cal N}_r} \frac{|T_{i,r}|}{48}\biggl(-6Cv_{i}^2+4E_{i}-6Cv_{i+1}^2+4F_{i} \biggl )\\
\ds \ge  \int_{\cup_{i \in {\cal N}_r} T_{i,r}} u^2 +  { 1 \over 48}\sum_{i \in {\cal N}_r} (|T_{i,r}|(-6Cv_{i}^2+4E_{i})+|T_{i-1,r}|(-6Cv_{i}^2+4F_{i-1})).
  \end{array}
\end{equation}
The condition~\eqref{eq:ABbraess} leads to the inequality 
\[
|T_{i,r}|E_i+|T_{i-1,r}|F_{i-1} = (|T_{i,r}|+|T_{i-1,r}|)v_i^2
\]
(see equation after (3.19) in \cite{BraessDahmen}), so that we get:
\bea
|T_{i,r}|(-6Cv_{i}^2+4E_{i})+|T_{i-1,r}|(-6Cv_{i}^2+4F_{i-1})
=(|T_{i,r}|+|T_{i-1,r}|)(4-6C)v_{i}^2. \nonumber
\eea
This term cancels for $C=2/3$. Hence for $0<C\le 2/3$, inequality \eqref{eq:sumqtir} becomes:
\be
\sum_{i \in {\cal N}_r} Q_{T_{i,r}}\ge \int_{\cup_{i \in {\cal N}_r} T_{i,r}}u^2. \nonumber
\ee
Therefore, for $0<C\le 2/3$ and $0<c\le 1$, 
\be
\sum_{i \in {\cal N}_r} Q_{T_{i,r}}\ge c\int_{\cup_{i \in {\cal N}_r} T_{i,r}}u^2. \nonumber
\ee
\paragraph{Triangles having all three vertices on the boundary}
Let $T_{c'}$ be such a triangle (see figure~\ref{fig:triangles}). We have to control:
\[
Q_{T_{c'}}=-\frac{C}{2} |T_{c'}| |u_{i+1}|^2
\]
by the integrals over the two triangles $T_{i,c}$ and $T_{j',c}$ surrounding $T_{c'}$. This can be achieved using the assumption 
\[
\frac{1}{24}\min(|T_{i,c}|,|T_{j',c}|) > \frac{C}{2}|T_{c'}|
\]
and using that from \eqref{eq:Qtic}, we have 
\[
Q_{T_{i,c}\cup T_{j',c}} \ge \min(|T_{i,c}|,|T_{j',c}|) \frac{u_{i+1}^2}{12}.
\]

In conclusion, we have that (\ref{eq:stabesti}) holds with $c=1/4$ for $0< C\le 2/3$. 

\section{Numerical results}
On the unit square $\Omega=(0,1) \times (0,1)$ we consider the problem
\begin{eqnarray*}
  (Id - \Delta)u(x,y) &=& x^3(y^2-2)-6xy^2+(1+x^2+y^2)sin(xy), 
  \quad (x,y) \in \Omega,\\
   u &=& x^3y^2+sin(xy), \quad (x,y) \in \partial{\Omega},
\end{eqnarray*}
whose exact solution is $u(x,y)=x^3y^2+sin(xy)$. We decompose
the unit square into non-overlapping subdomains 
with meshes generated in an independent manner. The computed solution is the solution at convergence of the discrete
algorithm (\ref{algo_discret})-(\ref{CI_discret}), with a
stopping criterion on the jumps of interface conditions that must be
smaller than $10^{-8}$. 
\\
\begin{rem}
In the implementation of the method, the main difficulty lies in computing
projections between non matching grids. In \cite{GJMN} we present 
an efficient algorithm in two dimensions to perform the
required projections between arbitrary grids, in the same spirit as in \cite{Gander:2001:OSW} for finite volume discretization with projections on
piecewise constant functions.
\end{rem}
%
\subsection{Choice of the Robin parameter $\alpha$}
In our simulations
the Robin parameter is either an arbitrary constant or is obtained by minimizing
the convergence rate (and depend of the mesh size in that case).
In the conforming two subdomains case, with constant mesh size $h$, the optimal theoretical value of
$\alpha$ which minimizes the
convergence rate at the continuous level is :
$$\alpha_{opt}=[(\pi^2+1)((\frac{\pi}{h})^2+1)]^{\frac{1}{4}}.$$
In the non-conforming case, has the mesh size is different for each side
of the interface, we consider the following values :
$$\alpha_{min}=[(\pi^2+1)((\frac{\pi}{h_{min}})^2+1)]^{\frac{1}{4}}$$
$$\alpha_{mean}=[(\pi^2+1)((\frac{\pi}{h_{mean}})^2+1)]^{\frac{1}{4}}$$
$$\alpha_{max}=[(\pi^2+1)((\frac{\pi}{h_{max}})^2+1)]^{\frac{1}{4}}$$
where $h_{min}$, $h_{mean}$ and $h_{max}$ stands respectively for the smallest
meanest or highest step size on the interface.
%
\subsection{An example of computed solution}
\begin{figure}[H]
  \centering 
  \includegraphics[height=6cm]{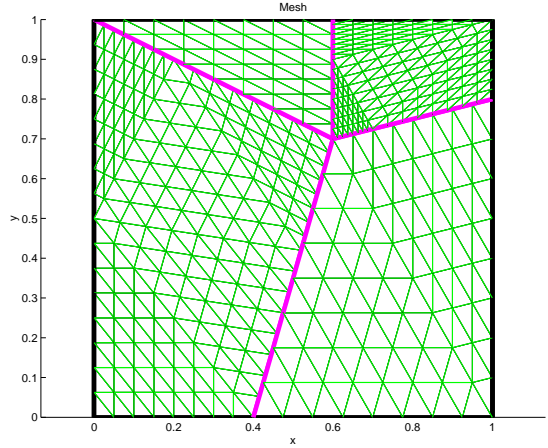}
  \caption{Domain decomposition with non-conforming meshes.}
  \label{fig:1}
\end{figure}

\begin{figure}[H]
  \centering 
  \includegraphics[height=8cm]{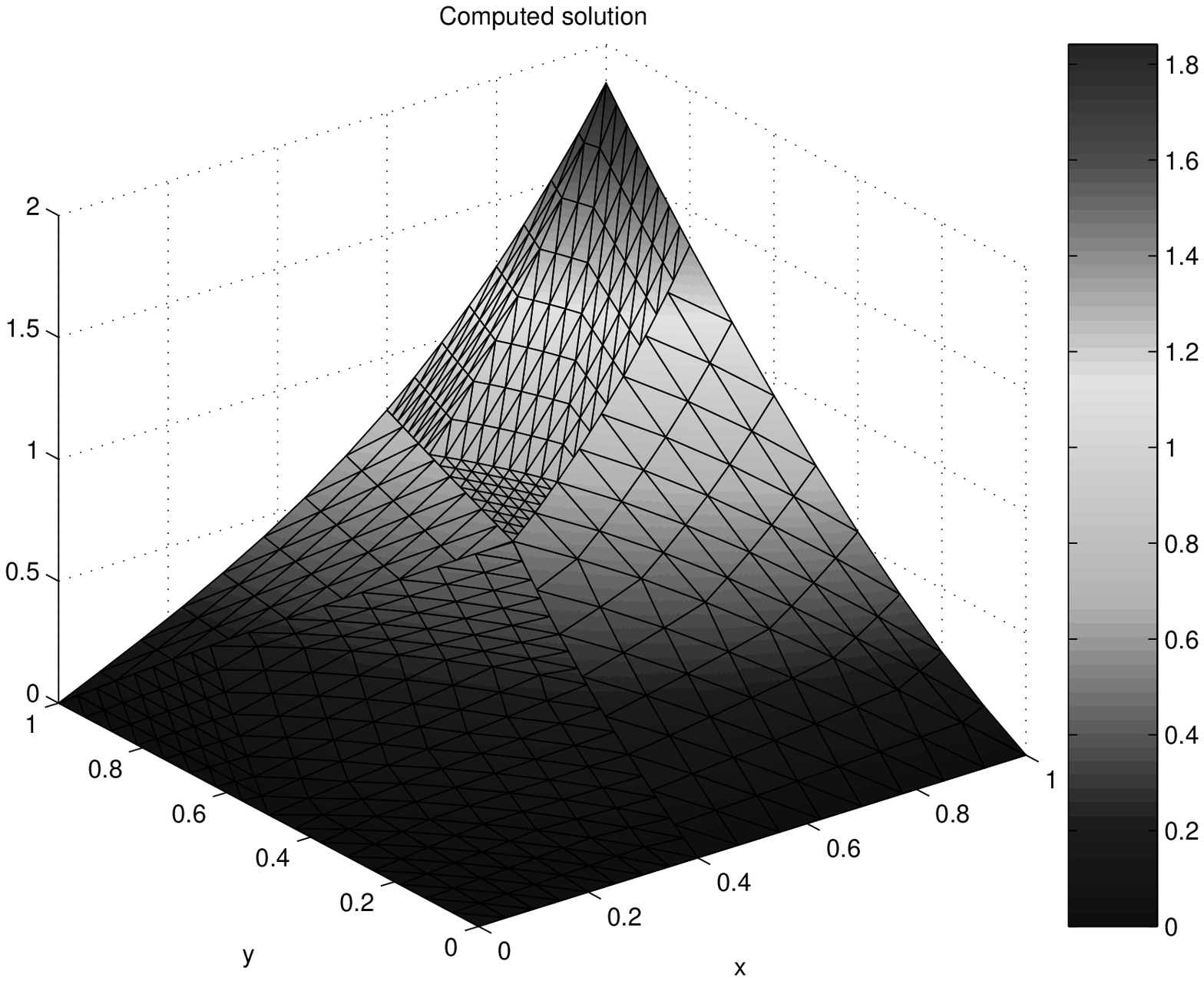}
  \caption{Computed solution.}
  \label{fig:2}
\end{figure}

%
We decompose the unit square into four non-overlapping subdomains 
with meshes generated as shown in Figure \ref{fig:1}.
The Robin parameter is $\alpha=10$.
On Figure \ref{fig:2} we show that the
computed solution is close to the continuous solution.
%
%
\subsection{$H^1$ error between the continuous and discrete solutions}
\begin{figure}[H]
  \centering 
  \includegraphics[height=4cm]{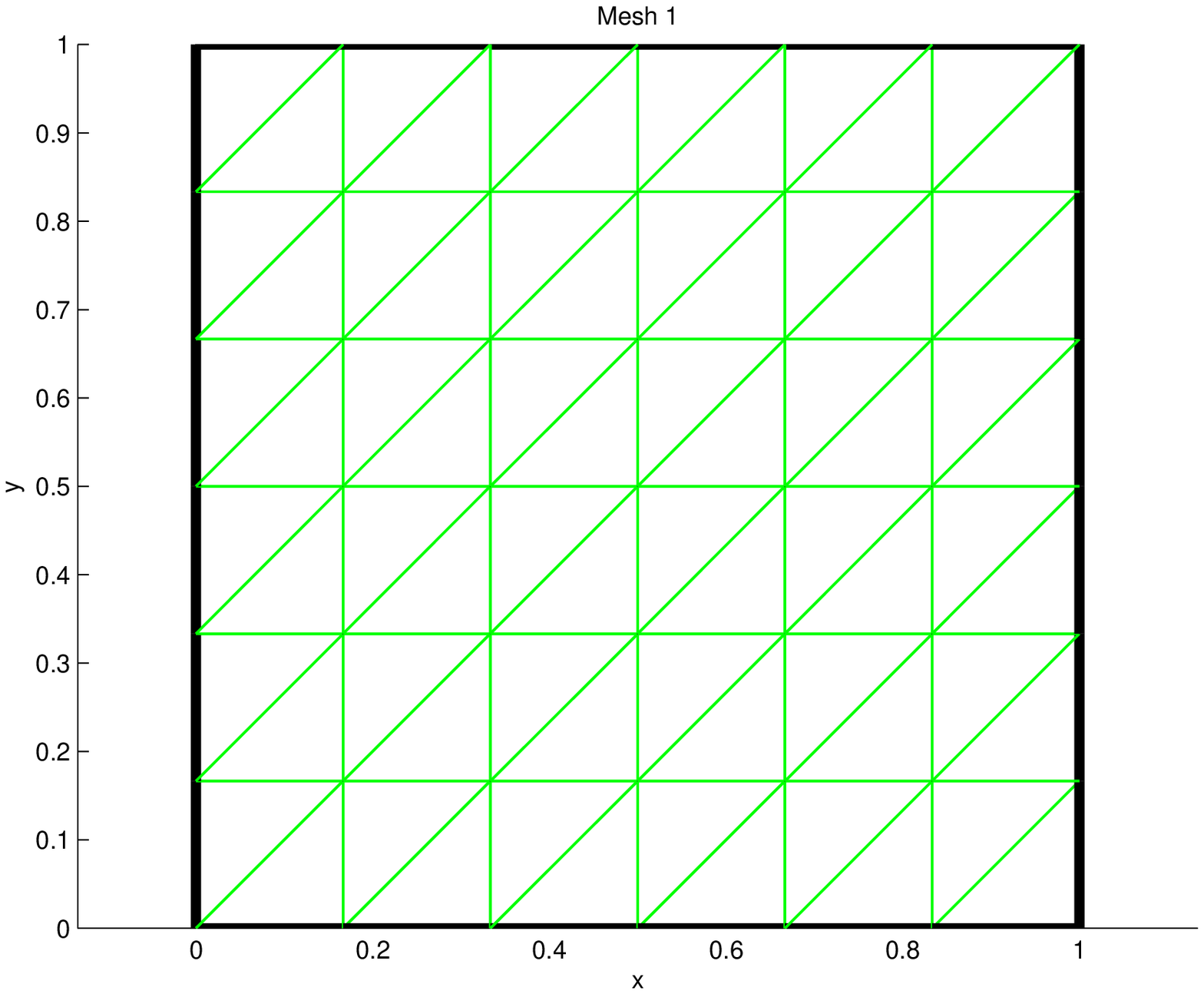}
  \includegraphics[height=4cm]{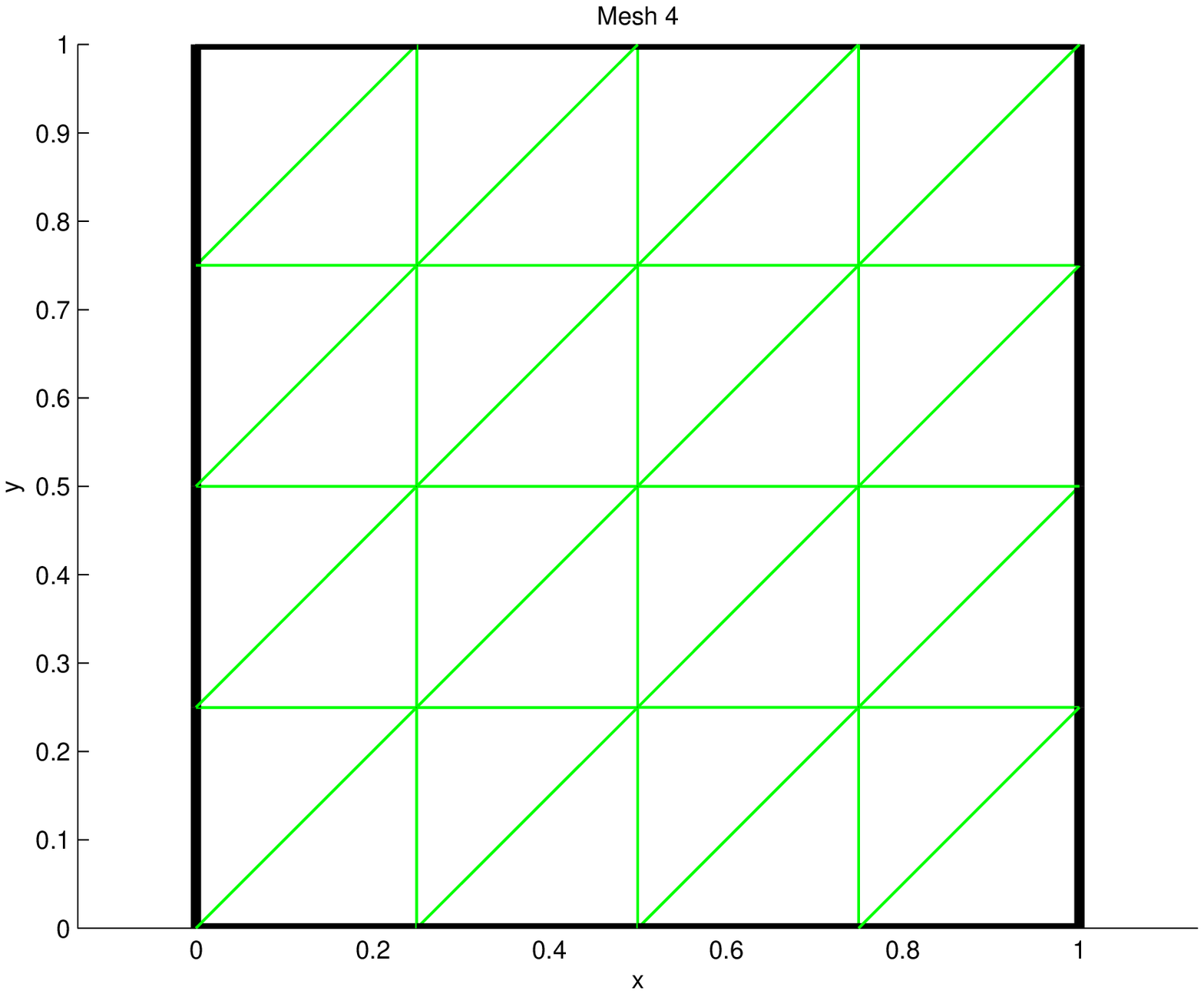}
  \caption{Uniform conforming meshes : mesh 1 (on the left), and mesh 4 (on the right)}
  \label{fig:3}
\end{figure}

\begin{figure}[H]
  \centering 
  \includegraphics[height=4cm]{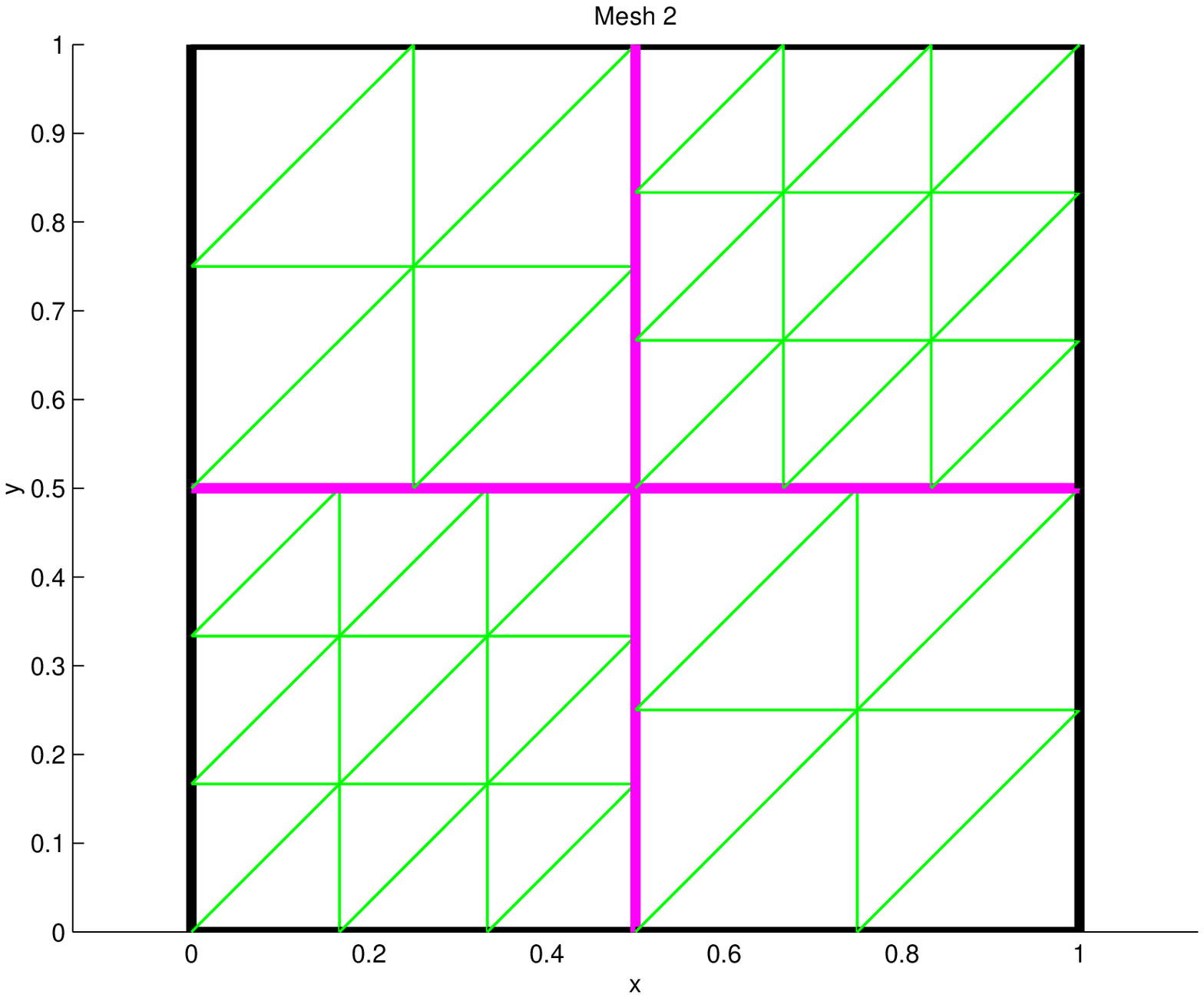}
  \includegraphics[height=4cm]{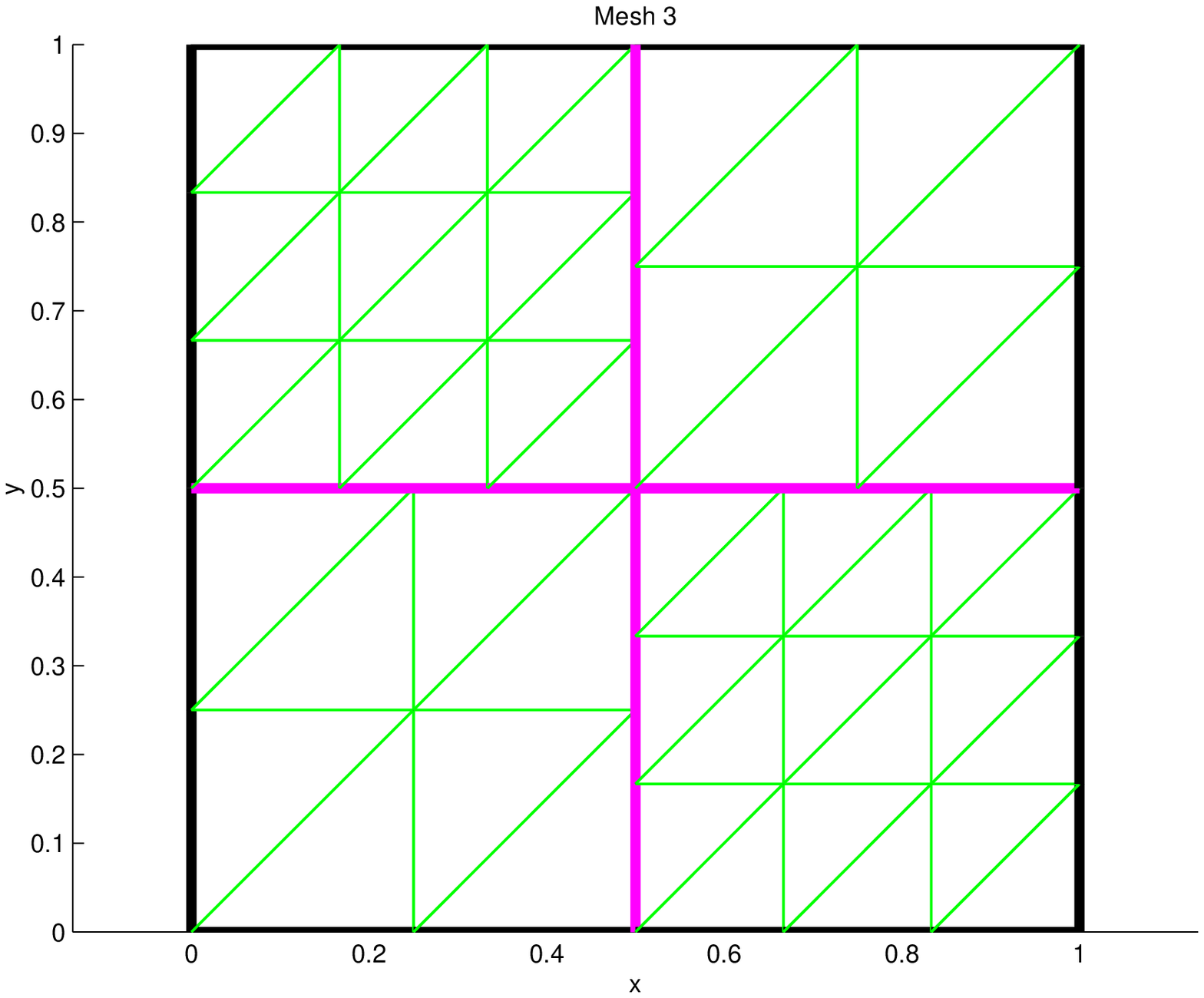}
  \caption{Non-conforming meshes : mesh 2 (on the left), and mesh 3 (on the right)}
  \label{fig:4}
\end{figure}

\begin{figure}[H]
  \centering 
  \psfrag{O1}{$\displaystyle\Omega^1$} 
  \psfrag{O2}{$\displaystyle\Omega^2$}
  \psfrag{O3}{$\displaystyle\Omega^3$}
  \psfrag{O4}{$\displaystyle\Omega^4$}
  \includegraphics[height=3cm]{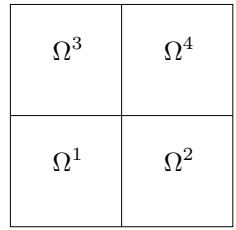}
  \caption{Non-overlapping domain decomposition of the unit square}
  \label{fig:5}
\end{figure}

\begin{figure}[H]
  \centering 
  \includegraphics[height=7cm]{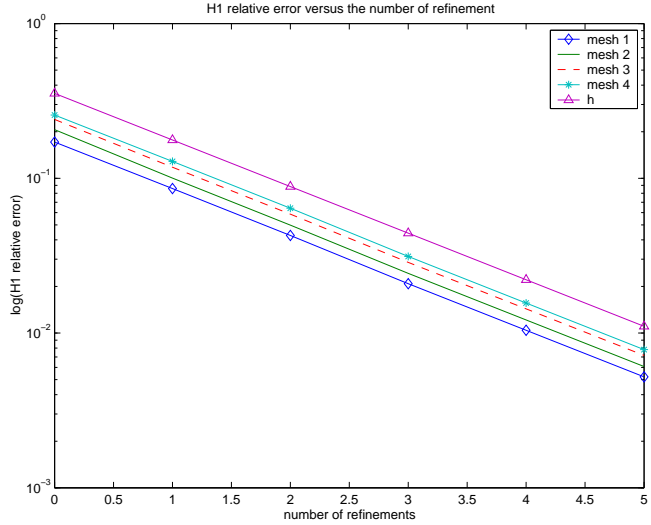}
  \caption{relative $H^1$ error versus the number of refinements for the initial meshes :
mesh 1, (diamond line), mesh 2 (solid line), mesh 3 (dashed line), and mesh 4 
(star line). The triangle line is the mesh size $h$ versus the number of refinements, in logarithmic scale}
  \label{fig:6}
\end{figure}

%
%
In this part, we compare the relative $H^1$ error in the non-conforming case
to the error obtained on a uniform conforming grid.\\\\
{\bf Definition of the relative $H^1$ error} : 
Let $K$ be
the number of subdomains. Let
$u_i=u_{|\Omega^i}, \ 1 \le i \le K$ (where u is the continuous
solution), and let $(\uu_h)_i=(\uu_h)_{|\Omega^i}$ where $\uu_h$ is the 
solution of the discrete problem (\ref{pbdiscret}). 
Now, let $E_{ex}=\|u\|_*$ and let $E_i=\|(\uu_h)_i-u_i\|_{H^1(\Omega^i)}, \ 1 \le i \le K$. 
Let $$E=(\sum_{i=1}^K E_i^2)^{1/2}.$$
The relative $H^1$ error is then $E/E_{ex}$.\\
In this example, we take $\alpha=\alpha_{mean}$ for the Robin parameter. This choice
is motivated by the results of section \ref{sec.Robinparam}, but we obtain similar results 
in the case of $\alpha=10$.
\\
We consider four initial meshes :
the two uniform conforming meshes (mesh 1 and 4) of figure \ref{fig:3}, and the two non-conforming meshes (mesh 2 and 3) of
figure \ref{fig:4}.
In the non-conforming case, the unit square is decomposed into 
four non-overlapping subdomains numbered as in figure \ref{fig:5}.
%
Figure \ref{fig:6} shows the relative $H^1$ error versus the number of
refinement for these four meshes, and the mesh size $h$ versus the number of
refinement, in logarithmic scale. At each refinement, the mesh size is
divided by two. The results of figure \ref{fig:6} show that the relative $H^1$ 
 error tends to zero at the same rate than the mesh size, and this
fits with the theoretical error estimates of theorem \ref{best-fit}.
On the other hand, we observe that the two curves corresponding to the 
non-conforming meshes (mesh 2 and mesh 3) are between the curves 
of the conforming meshes (mesh 1 and mesh 4). The relative $H^1$ error
for mesh 2 is smaller than the one corresponding to mesh 3, and this
is because mesh 2 is more refined than mesh 3 in subdomain $\Omega^4$, 
where the solution steeply varies.
%
%
More precisely, let us compare
for mesh 2, the relative $H^1$ error in the domain $\Omega^1\cup \Omega^2\cup\Omega^3$ to the relative $H^1$ error in 
the subdomain $\Omega^4$ (which is the subdomain where the solution steeply 
varies). This comparison can be done in Table 1. 
\medskip
\begin{center}
\begin{tabular}{|c|c|c|c|}
\hline
%
Refinement & $(E_1^2+E_2^2+E_3^2)^{1/2}/E_{ex}$ & $E_4/E_{ex}$ & $E/E_{ex}$\\
\hline
 0 & 1.45e-01 & 1.46e-01 & 2.06e-01\\
\hline
 1 & 7.17e-02 & 7.02e-02 & 1.004e-01\\
\hline
 2 & 3.59e-02 & 3.49e-02 & 5.01e-02\\
\hline
 3 & 1.79e-02 & 1.73e-02 & 2.49e-02\\
\hline
 4 & 8.73e-03 & 8.46e-03 & 1.21e-02\\
\hline
\end{tabular}
\end{center}
\begin{center}
{Table 1: Comparison, in the case of {\bf mesh 2}, for different refinements (column one), of the relative $H^1$ error in the domain composed by subdomains
$\Omega^1, \ \Omega^2$ and $\Omega^3$ (column 2) to the relative $H^1$ error 
in the subdomain $\Omega^4$ (column 3). The fourth column is the relative $H^1$ error in the whole domain.}
\end{center}
We observe that,
as expected, the relative $H^1$ error in the domain composed by subdomains
$\Omega^1, \ \Omega^2$ and $\Omega^3$ (second column of table 1)
is close to the relative $H^1$ error in 
the subdomain $\Omega^4$ (third column of table 1). Indeed, the mesh 2
is more refined in the subdomain $\Omega^4$ where the solution steeply 
varies.
\medskip
\begin{center}
\begin{tabular}{|c|c|c|c|}
\hline
%
Refinement & $(E_1^2+E_2^2+E_3^2)^{1/2}/E_{ex}$ & $E_4/E_{ex}$ & $E/E_{ex}$\\
\hline
 0 & 1.26e-01 & 2.04e-01 & 2.40e-01\\
\hline
 1 & 5.57e-02 & 1.04e-01 & 1.18e-01\\
\hline
 2 & 2.74e-02 & 5.22e-02 & 5.90e-02\\
\hline
 3 & 1.36e-02 & 2.59e-02 & 2.93e-02\\
\hline
 4 & 6.64e-03 & 1.26e-02 & 1.43e-02\\
\hline
\end{tabular}
\end{center}
\begin{center}
{Table 2: Comparison, in the case of {\bf mesh 3}, for different refinements (column one), of the $H^1$ relative error in the domain composed by subdomains
$\Omega^1, \ \Omega^2$ and $\Omega^3$ (column 2) to the $H^1$ relative error 
in the subdomain $\Omega^4$ (column 3). The fourth column is the $H^1$ relative error in the whole domain.}
\end{center}
\medskip
Let us now do the same comparison in the case of mesh 3. This mesh is
coarser in the subdomain $\Omega^4$ where the solution steeply 
varies. In table 2, we observe that
as expected, the $H^1$ relative error in the domain composed by subdomains
$\Omega^1, \ \Omega^2$ and $\Omega^3$ (second column of table 2)
is smaller (almost half) than the $H^1$ relative error in 
the subdomain $\Omega^4$ (third column of table 2). That one is close
to the $H^1$ relative error in 
the whole domain (fourth column of table 2), because mesh 3 is
coarser in the subdomain $\Omega^4$ where the solution steeply varies.
%
%
%
\subsection{Convergence : Choice of the Robin parameter}\label{sec.Robinparam}
%
Let us now study the convergence speed to reach the discrete solution, for
different values of the Robin parameter $\alpha$.
We first consider a domain decomposition in two subdomains, and then
in four subdomains. 
\subsubsection{2 subdomain case}
\begin{figure}[H]
  \centering
  \includegraphics[height=9cm]{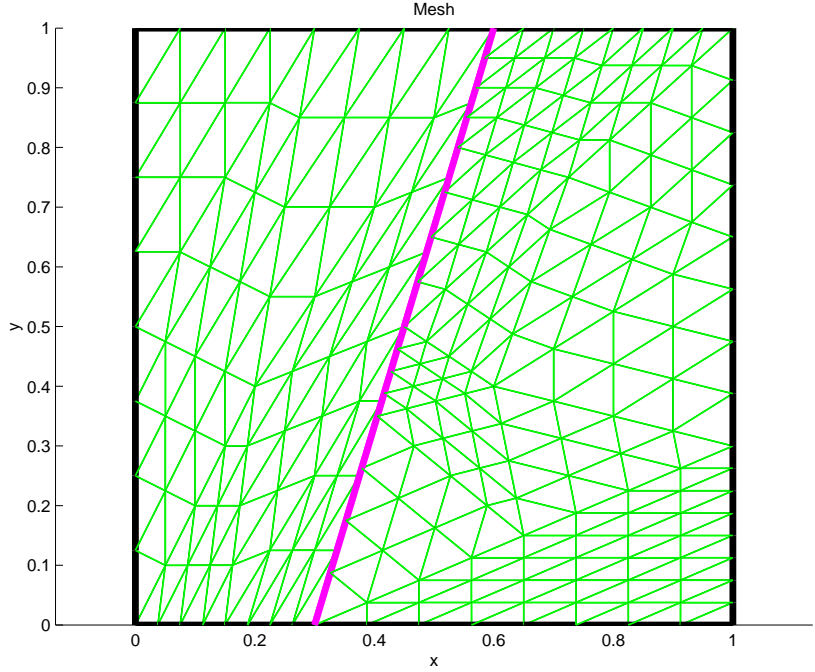}
  \caption{Domain decomposition in 2 subdomains with non-conforming grids}
\label{fig:7}       
\end{figure}

\begin{figure}[H]
  \centering
  \includegraphics[height=8cm]{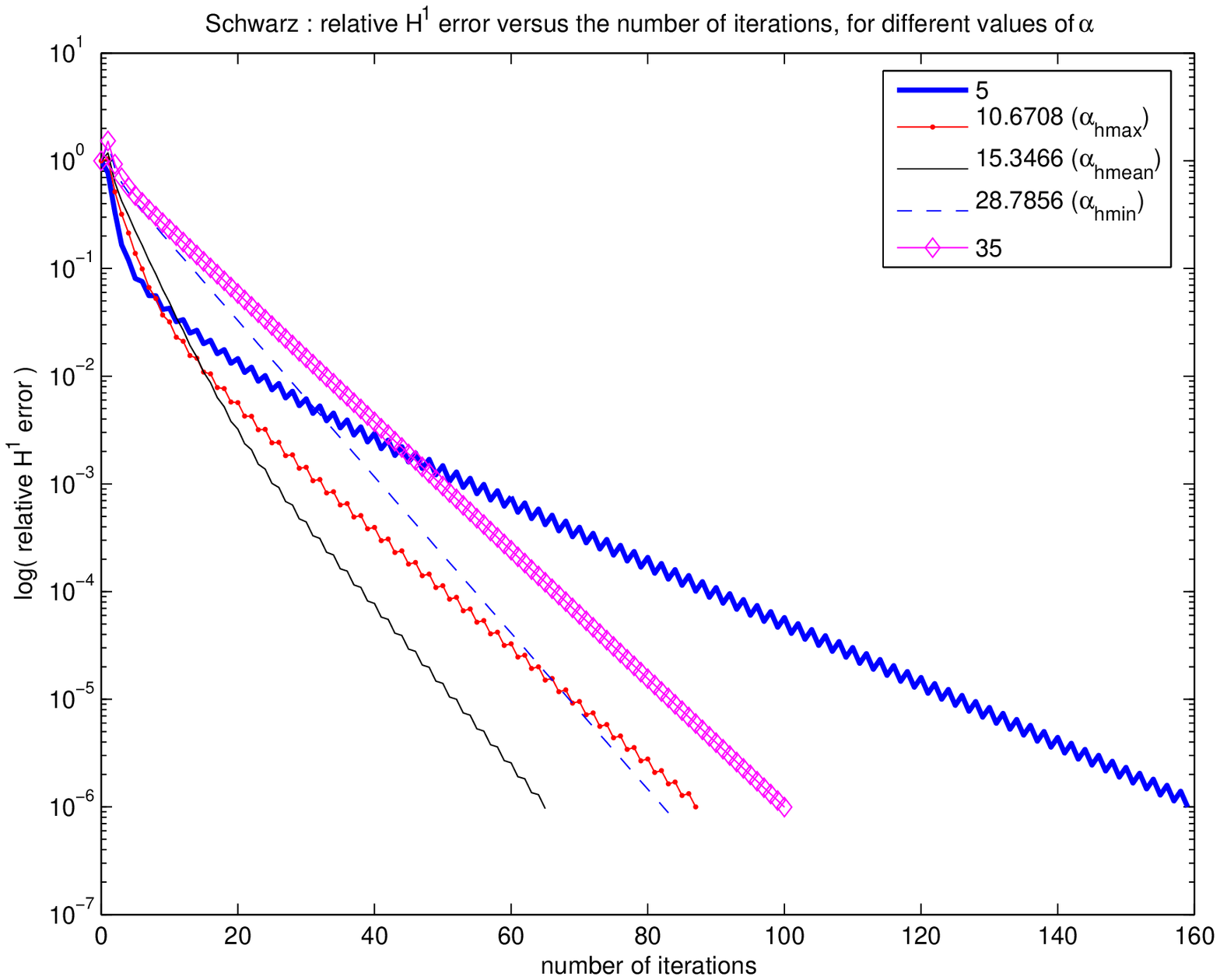}
  \caption{relative $H^1$ error between the discrete Schwarz converged 
solution and the iterate solution, for different values of the Robin 
parameter $\alpha$}
\label{fig:8}       
\end{figure}

\begin{figure}[H]
  \centering
  \includegraphics[height=8cm]{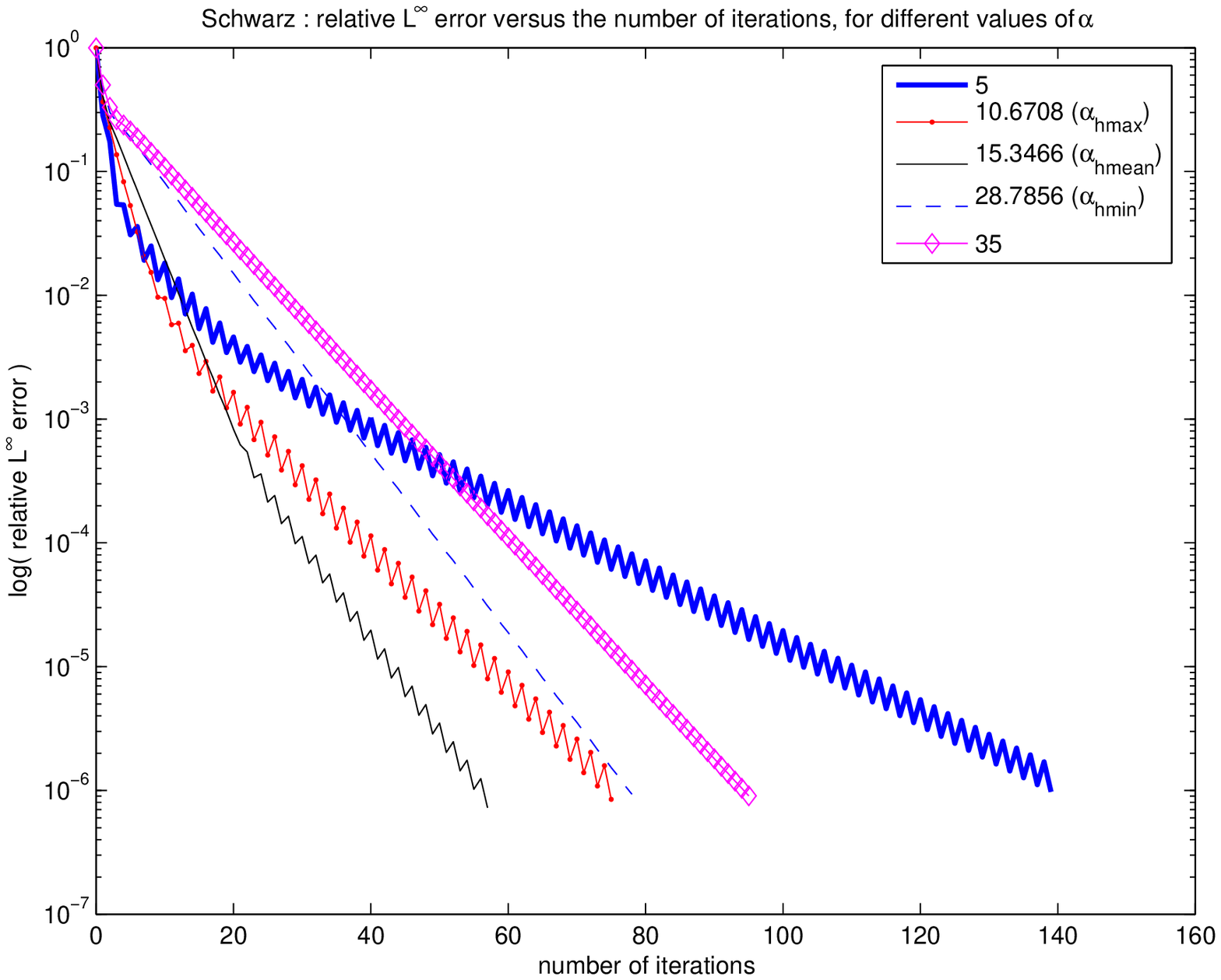}
  \caption{Relative $L^{\infty}$ error between the discrete Schwarz converged 
solution and
the iterate solution, for different values of the Robin parameter $\alpha$}
\label{fig:9}       
\end{figure}

\begin{figure}[H]
  \centering
  \includegraphics[height=8cm]{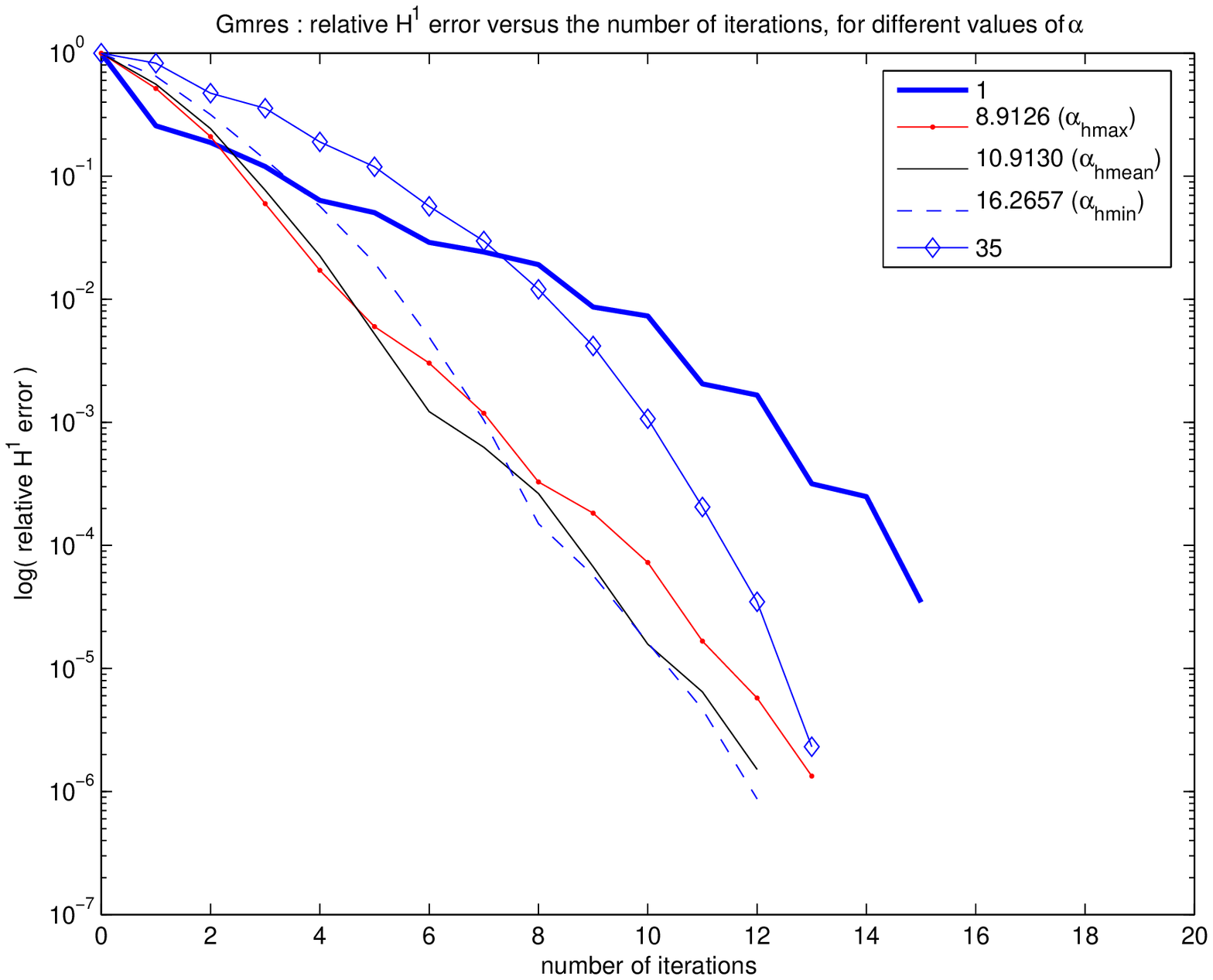}
  \caption{relative $H^1$ error between the discrete Gmres converged 
solution and the iterate solution, for different values of the Robin 
parameter $\alpha$}
\label{fig:10}       
\end{figure}

\begin{figure}[H]
  \centering
  \includegraphics[height=8cm]{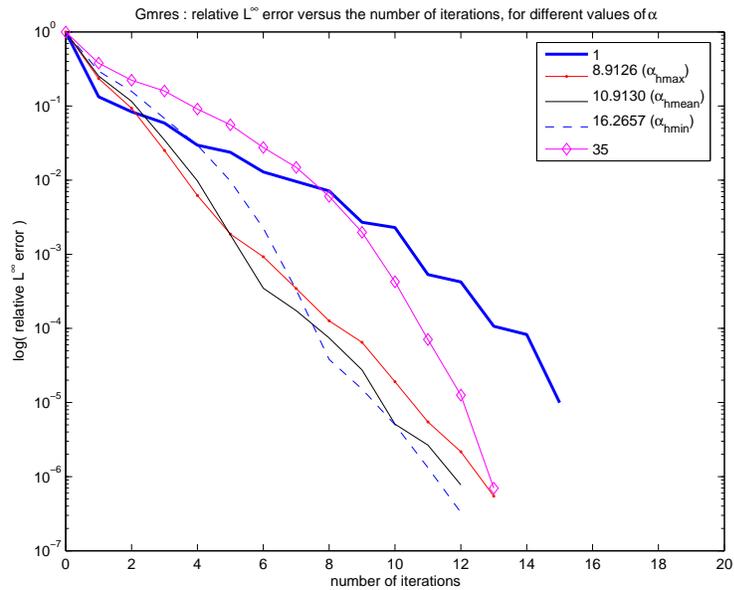}
  \caption{Relative $L^{\infty}$ error between the discrete Gmres converged 
solution and
the iterate solution, for different values of the Robin parameter $\alpha$}
\label{fig:11}       
\end{figure}

In this part, the unit square is decomposed in two subdomains
with non-conforming meshes (with $81$ and $153$ nodes respectively)
 as shown in figure \ref{fig:7}. On figure \ref{fig:8} we represent
the relative $H^1$ error between the discrete Schwarz converged 
solution and the iterate solution, for different values of the Robin 
parameter $\alpha$. We observe that
the optimal numerical value of the Robin parameter is close to $\alpha_{mean}$
and near $\alpha_{min}$ and $\alpha_{max}$.\\
\noindent
As the relative $H^1$ error didn't show where the error is highest,
we also look at the relative $L^{\infty}$ error between the discrete Schwarz 
converged solution and
the solution at iteration $p$, for different values of the Robin parameter $\alpha$.
We obtain similar results as for the relative $H^1$ error
(see figure \ref{fig:9}).\\
The Schwarz algorithm can be interpreted as a Jacobi algorithm applied
to an interface problem (see \cite{Nataf.4}). In order to accelerate the 
convergence, we can replace the Jacobi algorithm by a Gmres 
(\cite{Saad}) algorithm.
Figures \ref{fig:10} and \ref{fig:11} show respectively the
relative $H^1$ error and the relative $L^{\infty}$ error between the discrete 
Gmres converged solution and
the iterate solution, for different values of the Robin parameter $\alpha$.
In the case where $\alpha=\alpha_{mean}$, we observe that the convergence is 
accelerated by a factor 2 for Gmres, compared to Schwarz algorithm.
Also, the gap between the error values for different $\alpha$  
is decreasing when using Gmres algorithm, compared to Schwarz method.
The Gmres algorithm is less sensitive to the choice of the Robin
parameter.
\subsubsection{4 subdomain case}
\begin{figure}[H]
  \centering
  \includegraphics[height=9cm]{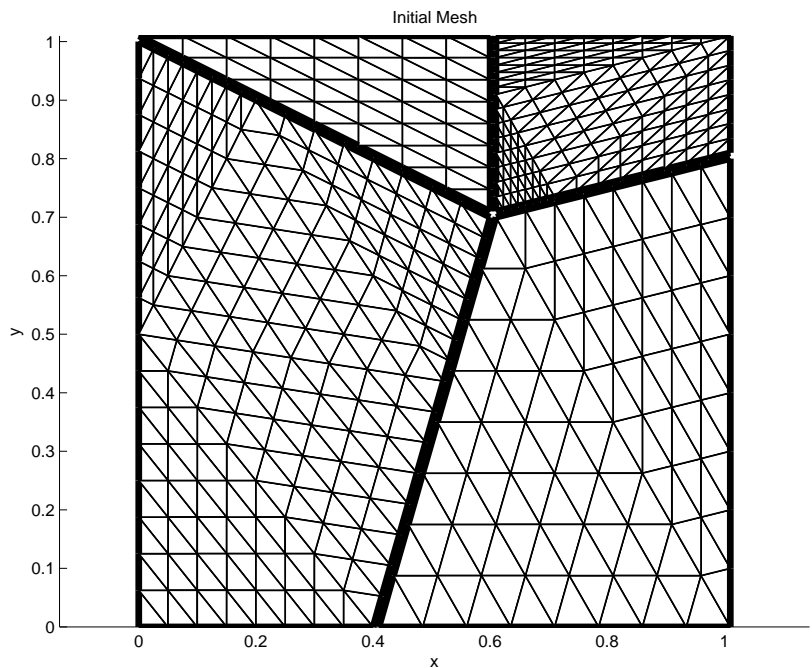}
  \caption{Domain decomposition in 4 subdomains with non-conforming grids}
\label{fig:12}       
\end{figure}

In this part, the unit square is decomposed into four subdomains
 with non-conforming 
meshes (with $189$, $81$, $45$ and $153$ nodes respectively) 
as shown in figure \ref{fig:12}. On figure \ref{fig:13} and 
\ref{fig:14} respectively, we represent
the relative $H^1$ error and the relative $L^{\infty}$ error between the discrete Schwarz converged 
solution and the iterate solution, for different values of the Robin 
parameter $\alpha$. We observe that
the optimal numerical value of the Robin parameter is close to $\alpha_{mean}$
and near $\alpha_{min}$ and $\alpha_{max}$, as in the two subdomain case.
%
\subsubsection{Conclusions}
%
The numerical results on the relative $H^1$ error between the continuous and 
discrete solutions correspond to
the theoretical error estimates of theorem \ref{best-fit}.
On the other hand, we observe that, for a fixed number of mesh points,
the relative $H^1$ error between the continuous and 
discrete solutions is smaller for a mesh refined in the region of the domain
where the solution steeply varies, than for a mesh which is coarser in that
region. In term of convergence speed to reach the discrete solution, the
Robin parameter $\alpha$ must depend of the mesh size, and our simulations show
that $\alpha=\alpha_{mean}$ is close to the optimal numerical value.
\begin{figure}[H]
  \centering
  \includegraphics[height=8cm]{errH1.eps}
  \caption{relative $H^1$ error between the discrete Schwarz converged 
solution and the iterate solution, for different values of the Robin 
parameter $\alpha$}
\label{fig:13}       
\end{figure}
\begin{figure}[H]
  \centering
  \includegraphics[height=8cm]{errLinf.eps}
  \caption{Relative $L^{\infty}$ error between the discrete Schwarz converged 
solution and
the iterate solution, for different values of the Robin parameter $\alpha$}
\label{fig:14}       
\end{figure}

\appendix

\section{Inf-sup condition.}
The purpose of this annex is to show that the proof of \cite{BB} can be
extended to the 3D situation. Indeed the main ingredients
required for the extensions have been proven in \cite{BraessDahmen}. 
Let us first recall a standard stability result in higher norms of the $L^2$
projection operator $\bar{\pi}_{k,\ell}$ from
$L^2(\Gamma^{k,\ell})$ onto
${\cal Y}_h^{k,\ell} \cap H^1_0(\Gamma^{k,\ell})$ orthogonal to $\tilde W_h^{k,\ell}$.
\begin{lem}\label{lemma:stability}
Making the hypothesis that the triangulation ${\cal T}_h^k$ is uniformly regular, there exists a constant $c>0$ such that
$$
\forall v \in H^{1 \over 2}_{00}(\Gamma^{k,\ell}), \ 
\|\bar{\pi}_{k,\ell}v\|_{H^{1 \over 2}_{00}(\Gamma^{k,\ell})} \le
c\|v\|_{H^{1 \over 2}_{00}(\Gamma^{k,\ell})}.
$$
\end{lem}
{\bf Proof} 
From (\ref{eq:stabesti}) we deduce a uniform inf-sup condition
between ${\cal Y}_h^{k,\ell} \cap H^1_0(\Gamma^{k,\ell})$ and 
$\tilde W_h^{k,\ell}$ in $L^2(\Gamma^{k,\ell})$.
It results that the projection operator $\bar{\pi}$ is stable in $L^2(\Gamma^{k,\ell})$
and thus there exists a constant $c_1>0$ such that
$$
\forall v \in H^{1 \over 2}_{00}(\Gamma^{k,\ell}), \ 
\|v-\bar{\pi}_{k,\ell}v\|_{L^2(\Gamma^{k,\ell})}
\le c_1 h^{1 \over 2}\|v\|_{H^{1 \over 2}_{00}(\Gamma^{k,\ell})}.
$$
Let $\tilde{\pi}_{k,\ell}$ denote the orthogonal projection operator from
$H^{1 \over 2}_{00}(\Gamma^{k,\ell})$ onto
${\cal Y}_h^{k,\ell} \cap H^1_0(\Gamma^{k,\ell})$ for
$H^{1 \over 2}_{00}(\Gamma^{k,\ell})$ inner product.
Then, for all $v$ in $H^{1 \over 2}_{00}(\Gamma^{k,\ell})$,
$$
\|\bar{\pi}_{k,\ell}v\|_{H^{1 \over 2}_{00}(\Gamma^{k,\ell})}
\le \|\tilde{\pi}_{k,\ell}v\|_{H^{1 \over 2}_{00}(\Gamma^{k,\ell})}
+\|\bar{\pi}_{k,\ell}v-\tilde{\pi}_{k,\ell}v\|_{H^{1 \over 2}_{00}(\Gamma^{k,\ell})}
.$$
Then, with an inverse inequality, there exists a constant $c_2>0$ such that
$$
\|\bar{\pi}_{k,\ell}v\|_{H^{1 \over 2}_{00}(\Gamma^{k,\ell})}
\le \|v\|_{H^{1 \over 2}_{00}(\Gamma^{k,\ell})}
+c_2 h^{-{1 \over 2}}\|\bar{\pi}_{k,\ell}v-\tilde{\pi}_{k,\ell}v\|_{L^2(\Gamma^{k,\ell})}
.$$
Thus,
$$
\|\bar{\pi}_{k,\ell}v\|_{H^{1 \over 2}_{00}(\Gamma^{k,\ell})}
\le \|v\|_{H^{1 \over 2}_{00}(\Gamma^{k,\ell})}
+c_2 h^{-{1 \over 2}}c'h^{1 \over 2}\|v\|_{H^{1 \over 2}_{00}(\Gamma^{k,\ell})},
$$
and then, with $c=1+c'c_2$ ,we have 
$$
\|\bar{\pi}_{k,\ell}v\|_{H^{1 \over 2}_{00}(\Gamma^{k,\ell})}
\le c\|v\|_{H^{1 \over 2}_{00}(\Gamma^{k,\ell})},
\forall v \in H^{1 \over 2}_{00}(\Gamma^{k,\ell}),
$$
which ends the proof of lemma
\ref{lemma:stability}.\\\\
Then from the definition of the $H^{-{1 \over 2}}_*(\Gamma^{k,\ell})$
norm, for any $p_{h,k,\ell}$ in $\tilde W_h^{k,\ell}$, there exists an element
$w^{k,\ell}$ in $H^{1 \over 2}_{00}(\Gamma^{k,\ell})$ such that
\bea
\int_{\Gamma^{k,\ell}} p_{h,k,\ell} w^{k,\ell} 
=_{(H^{1/2}_{00}) ^{\prime}(\Gamma^{k,\ell})}<p_{h,k,\ell},w^{k,\ell} >_{H^{1/2}_{00}(\Gamma^{k,\ell})}
=\|  p_{h,k,\ell}\|_{(H^{1 \over 2}_{00}(\Gamma^{k,\ell}))^{\prime}}\|w^{k,\ell} \|_{H^{1 \over 2}_{00}(\Gamma^{k,\ell})},\nonumber
\eea
and $w^{k,\ell} $ can be chosen such that
\bea
\|w^{k,\ell} \|_{H^{1 \over 2}_{00}(\Gamma^{k,\ell})}
=\|  p_{h,k,\ell}\|_{(H^{1 \over 2}_{00}(\Gamma^{k,\ell}))^{\prime}}.
\nonumber
\eea
We apply now the projection operator on $w^{k,\ell}$ from lemma \ref{lemma:stability}. We derive that 
$\bar{\pi}_{k,\ell}(w^{k,\ell})=w^{k,\ell}_h \in {\cal Y}_h^{k,\ell} \cap H^1_0(\Gamma^{k,\ell})$ and 
\bea
\|w^{k,\ell}_h \|_{H^{1 \over 2}_{00}(\Gamma^{k,\ell})}
\le c \|  p_{h,k,\ell}\|_{(H^{1 \over 2}_{00}(\Gamma^{k,\ell}))^{\prime}},
\nonumber
\eea
and
\bea
\int_{\Gamma^{k,\ell}} p_{h,k,\ell} w^{k,\ell}_h
= \int_{\Gamma^{k,\ell}} p_{h,k,\ell} w^{k,\ell}
=\|  p_{h,k,\ell}\|^2_{(H^{1 \over 2}_{00}(\Gamma^{k,\ell}))^{\prime}}.
\nonumber
\eea
It remains to lift $w^{k,\ell}_h$ over $\Omega^k$, this is done by
prolongating $w^{k,\ell}_h$ by zero over $\partial \Omega^k \setminus \Gamma^{k,\ell}$
and lifting this element of $H^{1 \over 2}(\partial \Omega^k)$ over
$\Omega^k$ as proposed in \cite{BG}.

\section{Extension in 2D for high order approximations}
\label{sec:Ext2D}
\begin{lem}\label{lem:etapsi_base}
Let $1 \le p\le 13$ be an integer. There exists $c$ and $C>0$ such that for all
$\eta\in\PP^p([-1,1])$ s.t.
$\eta(-1)=0$ there exists $\psi\in\PP^{p-1}([-1,1])$ s.t.
\[
\eta(1)=\psi(1)
\]
and
\[
J(\psi;\eta):=\int_{-1}^1 (\eta\,\psi -\frac{1}{4}(\eta-\psi)^2)
\ge  c \int_{-1}^1 \eta^2
\]
and
\[
\int_{-1}^1 \psi^2 \le C \int_{-1}^1 \eta^2 .
\]
\end{lem}
This lemma has been proven in the case $p=1$ at section \ref{sec.bestfit2D}.
For $p \ge 2$, we prove this lemma by studying for a given $\eta\in \PP^p([-1,1])$,
$\eta\neq 0$ the maximization problem : \\

Find $\psi\in\PP^{p-1}([-1,1])$ such that
\begin{equation}\label{eq:minconst}
J(\psi;\eta)=\max_{\begin{array}{l}
\phi \in\PP^{p-1}([-1,1])\\
\phi(1)=\eta(1)
\end{array}}J(\phi;\eta).
\end{equation}
The function $J$ is strictly concave and there exists a function
satisfying the constraint. This problem admits a solution. The functional
$J(\phi,\eta)$ being quadratic in $(\phi,\eta)$ and the constraint being
affine, the optimality condition shows that the problem reduces to a linear
problem whose right hand side depends linearly of $\eta$. The affine
constraint being of rank one, the problem (\ref{eq:minconst}) admits a
unique solution which depends linearly of $\eta$. Therefore, it makes
sense to introduce the operator:\\

\[
\begin{array}{rcl}
S: \PP^p_0([-1,1]) &\longrightarrow&  \PP^{p-1}([-1,1])\\
\hphantom{S: }\eta &\mapsto& \psi \mathrm{\ solution\ to\
(\ref{eq:minconst})}
\end{array}
\]
where $\PP^p_0([-1,1])$ is the set of functions of $\PP^p([-1,1])$ that
vanish at $-1$. In
Lemma~\ref{lem:etapsi_base}, we take $\psi=S(\eta)$.   The operator $S$ is
linear from a finite dimensional space to another so that it is continuous
for any norm on these spaces. Therefore there exists $C>0$ such that
$\int_{-1}^1
\psi^2 \le C \int_{-1}^1 \eta^2$.
Moreover, the function
\[
\begin{array}{rcl}
H: \PP^p_0([-1,1])\backslash \{0\} &\longrightarrow&  \R\\
\hphantom{S: }\eta &\mapsto& \ds\frac{J(S(\eta),\eta)}{\ds\int_{-1}^1
\eta^2}
\end{array}
\]
is continuous and such that $H(\eta)=H(\alpha\eta)$ for any $\alpha\neq
0$. Therefore, it reaches its minimum and proving
Lemma~\ref{lem:etapsi_base} amounts to prove
\begin{lem}\label{lem:etapsi}
Let $p\le 13$ and $\eta\in\PP^p([-1,1])$ s.t. $\eta(-1)=0$ and $\eta$ is not the null function.\\
Then,
\[
J(S(\eta);\eta) > 0.
\]
\end{lem}

{\bf Proof}.
We make use of the Legendre polynomials
\[
L_0(x)=1,\ L_1(x)=x,\ (m+1)L_{m+1}(x)=(2m+1)\,x\,L_m(x) - mL_{m-1}(x),\ m\ge 1
\]
Let us recall that for any $m\ge 0$,
\[\begin{array}{l}
L_m(1)=1,\ L_m(-1)=(-1)^m \\
\int_{-1}^1 L_m(x)\,L_{m'}(x)\,dx = \delta_{m\,m'}\ds\frac{2}{2m+1}
\end{array}\]
The polynomial $\eta$ is decomposed on the Legendre polynomials
\[
\eta=\sum_{m=1}^p \eta_m(L_m+L_{m-1})
\]
and $\psi=S(\eta)$ is sought in the form
\[
\psi=\sum_{m=0}^{p-1} \psi_m L_m
\]
so that it maximizes the quantity $J(\psi;\eta)$ under the constraint
$\eta(1)=\psi(1)$. This corresponds to the min-max problem
\[
\max_{\psi\in\PP^{p-1}([-1,1])}\min_{\mu\in\R} {\cal L}(\psi,\mu)
\]
where
\[
{\cal L}(\psi,\mu) = J(\psi;\eta)-\mu (\psi(1)-\eta(1)).
\]
We have to prove that the optimal value is positive. The optimality
relations w.r.t
$\psi$ give
\[
\frac{3}{2}(\eta_m+\eta_{m+1})-\frac{1}{2}\psi_m=\mu \frac{2m+1}{2},\
1\le m\le p-1
\]
and
\[
\frac{3}{2}\eta_1-\frac{1}{2}\psi_0=\frac{\mu}{2}.
\]
Therefore, we get
\begin{equation}\label{eq:psi}
\psi = 3\eta -3\eta_p\,L_p - \mu R_{p-1}
\end{equation}
where $\ds R_{p-1}=\sum_{m=0}^{p-1}(2m+1)L_m$ and
$\|R_{p-1}\|_{L^2(]-1,1[)}^2=2p^2$. Hence, the dual problem writes
\[
\min_{\mu\in\R} G(\mu;\eta)
\]
where
\[
G(\mu;\eta) := J(3\eta -3\eta_p\,L_p - \mu
R_{p-1};\eta)-\mu(\psi(1)-\eta(1))
\]
and $\psi$ satisfies (\ref{eq:psi}).
After some calculations, we get
\begin{equation}\label{eq:G}
G(\mu;\eta) = \frac{p^2}{2} \mu^2 -\mu (2\eta(1)-3\eta_p)
+(2\|\eta\|_{L^2(]-1,1[)}^2-\frac{9}{2}\frac{\eta_p^2}{2p+1}).
\end{equation}
The leading coefficient of $G(\mu;\eta)$ is positive so that proving
$\min_\mu G(\mu;\eta)$ is positive (and hence Lemma~\ref{lem:etapsi}) is
equivalent to prove
\begin{lem}
For $p\le 13$, the discriminant of (\ref{eq:G}):
\begin{equation}\label{eq:discr}
\Delta(\eta) := (2\eta(1)-3\eta_p)^2 +
p^2 (-4\|\eta\|_{L^2(]-1,1[)}^2 + 9\frac{\eta_p^2}{2p+1})
\end{equation}
  is negative if $\eta\in\PP^p([-1,1])$, $\eta(-1)=0$ and $\eta$ is not
the null function.
\end{lem}
{\bf Proof}.
We first treat separately the case $p=2$. In this case, a direct computation shows that 
\[
\Delta(\eta) =-\frac{80}{3} \eta_1^2-\frac{40}{3} \eta_2\eta_1-\frac{133}{15}\eta_2^2
\]
The discriminant of the corresponding bilinear form is $-8632/9$. It is negative and the lemma is proved in this case.\\
We consider now the case $p\ge 3$.
Let us introduce the vector space $Q^p=\{
\eta\in\PP^p([-1,1])
\hbox{ s.t. }\
\eta(-1)=0\}$.  The function $\Delta(\eta)$ is quadratic so that it
suffices to study the extrema of $\Delta(\eta)/\|\eta\|^2_{L^2(]-1,1[)}$
over $Q^p$ or equivalently to prove that the associated symmetric 
quadratic form in
negative, i.e. its eigenvalues are negative. They correspond to the Lagrange
multiplier solutions
$\mu_1$ of the following min-max problem
\begin{equation}\label{eq:minMaxDiscr}
\min_{\eta\in Q^p} \max_{\mu_1\in\R} {\cal L}_e(\eta,\mu_1)
\end{equation}
where
\[
{\cal L}_e(\eta,\mu_1) := \Delta(\eta) - \mu_1 (\|\eta\|_{L^2(]-1,1[)}^2-1).
\]
We have to prove that $\mu_1<0$. We have 
\[\begin{array}{l}
0=\ds<\frac{\partial{\cal L}_e}{\partial \eta},\delta \eta >\\
\ds\hphantom{ddd}=
2(2\eta(1)-3\eta_p)(2\delta\eta(1)-3\delta\eta_p)+p^2(-8<\eta,\delta\eta>
+18\frac{\eta_p\delta\eta_p}{2p+1})
-2\mu_1 <\eta,\delta\eta>
\end{array}\]
where $<\, ,\,>$ denotes the $L^2$ scalar product on $L^2(]-1,1[)$ and
$\delta\eta\in Q^p$.
\\ Let us consider the vector space $(1-x^2) \PP^{p-3}\subset Q^p$. Any
  function $\gamma$ in $(1-x^2) \PP^{p-3}$ satisfies $\gamma(-1)=\gamma(1)=0$
and $\gamma_p=0$.
The optimality relation w.r.t. to $(1-x^2) \PP^{p-3}$ gives
\[
(-8p^2-2\mu_1)<\eta,\delta\eta>=0,\ \ \forall
\delta\eta\in (1-x^2) \PP^{p-3}.
\]
We have either $\mu_1=-4p^2<0$ or  $\eta$  solution to (\ref{eq:minMaxDiscr})
belongs to the space $\{(1-x^2) \PP^{p-3}\}^\bot
\cap \PP^p$. The first case corresponds to a negative value for 
$\mu_1$ which is in
agreement with the lemma to be proved. Let us study the latter case. 
We shall make
use of
\begin{lem}\label{lem:LagrangePrime}
\begin{eqnarray}
\label{eq:LP1}
\int_{-1}^1 L'_m\,L'_{m'}\,(1-x^2) \, dx = 0, \ m\neq m',\\
\int_{-1}^1 {L'_m}^2= m(m+1),\\
\int_{-1}^1 L'_m\,L'_{m+1}= 0,\\
\int_{-1}^1 L'_{m-1}\,L'_{m+1}= m(m-1),\\
L'_m(-1) = (-1)^{m+1} \frac{m(m+1)}{2}.
\end{eqnarray}
see \cite{AS}. 
\end{lem}
 From Lemma~\ref{lem:LagrangePrime}, it can be proved that
\begin{lem}\label{lem:11}
\[
\{(1-x^2) \PP^{p-3}\}^\bot \cap \PP^p = \hbox{ Span}\{L_{p},\,L'_p,\,L'_{p-1}\}.
\]
\end{lem}
{\bf Proof.} From (\ref{eq:LP1}), it can be checked easily that
\[
\{(1-x^2) \PP^{p-3}\}^\bot \cap \PP^p = \hbox{ Span}\{L'_{p+1},\,L'_p,\,L'_{p-1}\}.
\]
Moreover, we have
\[
L'_{p+1}(x)=(2p+1)L_p(x)+L'_{p-1}(x))
\]
and thus lemma \ref{lem:11}.\\
Therefore, there exists $\lambda_1,\lambda_2,\lambda_3 \in\R$ s.t.
$\eta=\lambda_1 L_{p} + \lambda_2 L'_{p} + \lambda_3 L'_{p-1}$.
Since $\eta$ is defined up  to a constant and we only
have to consider the two cases
$\lambda_1=1$ or $\lambda_1=0$.\\
{\bf Case 1} $\lambda_1=1$\\
 From $\eta(-1)=0$, we get
\[
1-\lambda_2 \frac{p(p+1)}{2} +\lambda_3 \frac{p(p-1)}{2} =0,
\]
so that
\[
\lambda_2 = \frac{2}{p(p+1)} +\lambda_3  \frac{p-1}{p+1}.
\]

\begin{eqnarray*}
\Delta(\eta)=\lefteqn{  {\displaystyle - 4\,{\displaystyle \frac
{(p-1)\,p^{2}\,(p^{2} +  1)}{p + 1}}\lambda_2^{2}}  -  {\displaystyle
\frac {(24\,p^{4} - 20\,p^{3} - 8\,p^{2} + 4\,p)\,
  }{(p + 1)\,(2\,p + 1)}\lambda_2} } \\
  & & \mbox{} - {\displaystyle \frac {29\,p^{2} + 13\,p - 1 - p^{3
}}{(p + 1)\,(2\,p + 1)}} \mbox{\hspace{223pt}}
\end{eqnarray*}

Since $p$ is supposed larger than 1, the
leading coefficient of $\Delta(\eta)$ is negative. If the discriminant of
$\Delta(\eta)$ is negative, the polynomial is negative for any
$\lambda_2$. This discriminant has the value
\[
16\,{\displaystyle \frac {(p^{2} - 13\,p - 8)\,(p - 1)\,p^{3}}{2
\,p + 1}}
\]
and is negative for $2\le p \le 13$.

{\bf Case 2} $\lambda_1=0$\\
 From $\eta(-1)=0$, we get
\[
-\lambda_2 \frac{p(p+1)}{2} +\lambda_3 \frac{p(p-1)}{2} =0,
\]
so that
\[
\lambda_2 = \lambda_3  \frac{p-1}{p+1}.
\]
Since $\eta$ is an eigenvalue, it is not zero and the above relation
shows that we can take $\lambda_3=1$. Then, we have
$\lambda_2 = \frac{p-1}{p+1}$ so that
\[
\Delta(\eta)=\frac{-4(p-1)p^2(p^2+1)}{(p+1)} < 0.
\]
\vspace{2mm}

\indent {\bf Acknowledgments.} The authors would like to thank Martin J. Gander for his help in the implementation of the method, especially for computing
projections between arbitrary grids in two dimensions.

\end{document}